\documentclass[11pt]{article}

\usepackage[british]{babel}
\usepackage{float} 
\usepackage{tikz}
\usetikzlibrary{matrix,arrows,decorations.pathmorphing}

\usepackage{amsmath,amssymb,wasysym,amsfonts,amsthm,amsthm}
\usepackage{empheq}
\usepackage{bm}
\usepackage{ulem}

\usepackage{mathtools}
\mathtoolsset{showonlyrefs}

\usepackage{graphicx,epsfig}
\usepackage{subfig}
\usepackage{overpic}

\usepackage{txfonts}
\usepackage[T1]{fontenc}
\usepackage{textcomp}
\usepackage[utf8x]{inputenc}
\usepackage{ucs}
\usepackage[dvipsnames]{xcolor}
\usepackage{titletoc}
\usepackage{fancyhdr}
\usepackage[margin=10pt,font=small,textfont=it,labelfont=bf]{caption}
\usepackage{hyperref}
\usepackage{url,enumerate,framed,ragged2e}

\newtheorem{thm}{Theorem}
\newtheorem{lem}[thm]{Lemma}

\newtheorem{deftn}[thm]{Definition}
\newtheorem{ass}[thm]{Assumptions}
\newtheorem{res}[thm]{Result}
\theoremstyle{remark}
\newtheorem{rmk}[thm]{Remark}
\numberwithin{equation}{section}
\numberwithin{thm}{section}

\newcommand{\eps}{\varepsilon}

\definecolor{myyellow}{RGB}{217, 162, 79}
\definecolor{mypurple}{RGB}{115, 64, 120}
\definecolor{myblack}{RGB}{76, 108, 165}
\definecolor{mydarkblack}{RGB}{40, 52, 139}
\definecolor{myorange}{RGB}{244, 148, 15}
\definecolor{mydarkred}{RGB}{190, 22, 34}

\title{Travelling front solutions in a spatially heterogeneous reaction-diffusion system}
\author{M. Chirilus-Bruckner$^*$, L. van Vianen$^*$, F. Veerman\footnote{Mathematical Institute, Leiden University, Einsteinweg 55, 2333 CC, Leiden, the Netherlands}}
\date{}

\begin{document}

\maketitle

\vspace{-1cm}

\begin{abstract}
We investigate a two-component reaction-diffusion system with a slow-fast structure and spatially varying coefficients $f_1$ and $f_2$ appearing in the slow equation. Under mild boundedness and regularity conditions on $f_1$ and $f_2$ the system is shown to exhibit bi-stability in the form of two stable stationary heterogeneous background states. These background states can be connected by stationary and travelling front solutions. Travelling fronts feature an interface that moves with a non-uniform speed through the motionless spatially varying background states it connects.  We construct both the background states and stationary fronts using an extension of Fenichel theory to the non-compact case. Additionally, we establish the existence of travelling front solutions and derive a leading-order expression for the dynamic position of the moving interface through a time-dependent spatial dynamics approach. This expression takes the form of a delay-differential equation, and its accuracy is validated through numerical simulations. A key contribution of our work lies in the general treatment of $f_1$ and $f_2$, which are neither (necessarily) asymptotically small nor restricted to specific forms such as periodic or localized structures. Furthermore, our derivation of the front position formula circumvents the traditional reliance on spectral analysis, enabling us to describe front dynamics beyond bifurcations from stationary fronts. This approach has the potential to be extended to other settings in which spectral properties at onset preclude conventional reduction techniques.

\end{abstract}

{\footnotesize
\tableofcontents
}

\section{Introduction}
We consider travelling front solutions of the reaction-diffusion equation
\begin{align}\label{eq:PDE_model}
\left\{
\begin{aligned}
\partial_t U &= \varepsilon^2\partial_x^2 U + U - U^3 -  \varepsilon(\alpha V +\gamma) \, ,\\[.2cm]
\tau \partial_t V &= \partial_x^2 V-[1+f_1(x)]V +[1+ f_2(x)]U \, ,
\end{aligned}
\right.
\end{align}
with parameters $\alpha, \gamma \in \mathbb{R}$, $\tau> 0$, perturbation parameter $0<\varepsilon \ll 1$ and spatially varying coefficients $f_1$ are $f_2$ that are not necessarily assumed to be small or with a specific structure (e.g. periodic or localized). 

\begin{ass}\label{ass:f} 
Let the coefficients $f_k, k = 1, 2$ in \eqref{eq:PDE_model} fulfil the following:
\begin{itemize}
        \item $f_k \in C^3_b(\mathbb{R})$, i.e. $f_k$ is bounded and sufficiently smooth, and
        \item $\inf_{x \in \mathbb{R}} [f_k(x) + 1] > 0$.
    \end{itemize}
\end{ass}
    
The present work establishes that \eqref{eq:PDE_model} exhibits bi-stability, characterized by two stable stationary background states that vary in space due to the presence of $f_1, f_2$ and further give rise to stationary and travelling front solutions that connect them (see Figure~\ref{fig:front_profile_numerics}).

\paragraph{Background: Nagumo/Allen-Cahn setting.} Consider the Nagumo equation
\begin{align*}
\partial_t U = \varepsilon^2\partial_x^2 U + U - U^3 { -} \varepsilon \gamma \, .
\end{align*}
Provided that $\varepsilon$ and $\gamma$ are such that $u-u^3+ \varepsilon \gamma=0$ has three distinct real roots, it has three constant background states $u^-(\varepsilon)< u^0(\varepsilon)< u^+(\varepsilon)$ with $u^\pm(\varepsilon)= \pm 1 - \frac{\gamma}{2}\varepsilon +\mathcal{O}(\varepsilon^2)$,  $u^0(\varepsilon)= \gamma \varepsilon + \mathcal{O}(\varepsilon^2)$. Furthermore, it
features uniformly moving front solutions
\begin{align*}
U_{TF}(x,t) = u_{*,+}(\varepsilon) + u_{*,-}(\varepsilon) \tanh\left[ \sqrt{2} u_{*,-}(\varepsilon) \left(\frac{x-\varepsilon^2c(\varepsilon)t}{2\varepsilon} \right)\right] \, ,
\end{align*}
with $u_{*,\pm}(\varepsilon) = (u^+(\varepsilon) \pm u^-(\varepsilon))/2$ whose speed is given by $c(\varepsilon)=- (u^+(\varepsilon)-2u^0(\varepsilon)+ u^-(\varepsilon ))/(\sqrt{2}\varepsilon) = \frac{3}{\sqrt{2}} \gamma + \mathcal{O}(\varepsilon)$. In particular, if $\gamma = 0$ there are only stationary front solutions and the forcing $\varepsilon \gamma$ induces the movement. 

\paragraph{Background: Constant-coefficient reaction-diffusion model.} A similar structure is inherited by the constant-coefficient version of \eqref{eq:PDE_model} 
\begin{align}\label{eq:PDE_model_const}
\left\{
\begin{aligned}
\partial_t U &= \varepsilon^2\partial_x^2 U + U - U^3 -   \varepsilon(\alpha V +\gamma) \, ,\\[.2cm]
\tau \partial_t V &= \partial_x^2 V-V +U \, .
\end{aligned}
\right.
\end{align}
It has stable background states $(u,v)^{\pm}(\varepsilon) = \pm (1,1) + \mathcal{O}(\varepsilon)$ that can be connected by stable stationary, uniformly travelling and dynamically travelling front solutions (cf. \cite{C-B.D.vH.R.2015,C-B.vH.I.R.2019,C-BvHR.2024}). More specifically, travelling front solutions $(U,V)(x,t) = (u_{\mathrm{tf}}, { u_{\mathrm{tf}}})(x-\varepsilon^2 ct)$ exist if the leading order existence condition
\begin{align}\label{eq:existence_const_coeff}
   \alpha v_* + \gamma  = \frac{\sqrt{2}}{3} c \, , \qquad v_* := \frac{c\hat{\tau}}{\sqrt{c^2 \hat{\tau}^2+4}} \, ,
\end{align}
relating the speed $c \in \mathbb{R}$ and the system parameters is fulfilled. Notably, we assumed a "large" relaxation parameter $\tau =: \frac{\hat{\tau}}{\varepsilon^2} > 0$. This shows, in particular, that if $\tau$ is not large enough, so $\hat{\tau} = 0$, then the second component does not impact the front motion. Moreover, from \cite{Peter2009,Peter2008,Peter2010} multi-front and pulse solutions of \eqref{eq:PDE_model_const} are known to have particularly rich dynamics in this regime. 

\paragraph{Fast-reaction scaling.} Given the role of $\tau = \frac{\hat{\tau}}{\varepsilon^2} > 0$ for the constant-coefficient setting we choose to stick to this parameter choice throughout this article which, in fact, can be related via $s = \varepsilon^2 t$ to the fast-reaction limit (\cite{KS.2025,STT.2024,TT.2024})
\begin{align}\label{eq:PDE_model_fast_reaction}
\left\{
\begin{aligned}
\varepsilon^2 \partial_s U &= \varepsilon^2\partial_x^2 U + U - U^3 -  \varepsilon(\alpha V +\gamma) \, ,\\[.2cm]
\hat{\tau} \partial_s V &= \partial_x^2 V-[1+f_1(x)]V +[1+ f_2(x)]U \, .
\end{aligned}
\right.
\end{align}

{ 
\subsection{Main results and Plan of the paper}
The first main contribution is to prove that, under the mild Assumptions~\ref{ass:f}, the PDE admits two stable stationary heterogeneous background states $(u_b^\pm(x;\varepsilon),v_b^\pm(x;\varepsilon))$. Their rigorous construction is sketched in Section~\ref{s:background_states}. Adopting a spatial-dynamics viewpoint, we analyse the steady-state ODE, augmented by the additional variable $\chi_x=1$, which features critical and slow manifolds that are non-compact. We sketch an extension of classical Fenichel theory to the non-compact setting and prove the existence of two bounded solutions that correspond to the two heterogeneous background states in the PDE, which are shown in \cite{LvV_thesis} to be stable.\\

\medskip

The second main contribution, given in Section~\ref{s:front_solutions}, is existence theory for front solutions connecting these motionless backgrounds. Stationary fronts are obtained in Section~\ref{ss:stationaryfronts} as heteroclinic connections between the two bounded solutions of the steady-state ODE from Section~\ref{s:background_states}, using a Melnikov-type computation that yields a leading-order existence and selection condition fixing admissible front locations, reflecting the loss of translation invariance caused by heterogeneity. Building on this, we further establish in Section~\ref{ss:travelingfronts} the existence of genuinely travelling fronts whose interface moves with non-uniform speed through stationary heterogeneous tails given by the heterogeneous background states from Section~\ref{s:background_states}. We sketch the rather lengthy rigorous construction from \cite{LvV_thesis}, in which one first controls the initial-value problem with carefully defined ``initialised fronts'' and then tracks the evolution of the time derivatives $(\partial_s U,\partial_s V)$, which decay to zero in the far field since the backgrounds do not move. We round off our exposition with Section~\ref{s:delay_equation}, where we derive a reduced law for the motion of the interface of travelling front solutions. In the singular limit $\varepsilon\to0$, the front position $z(s)$ satisfies a delay-differential equation in which the instantaneous drift $z'(s)$ depends on a history-dependent functional $W[z]$ arising from solving an auxiliary linear parabolic problem for the slow field. Crucially, this reduction is obtained without spectral analysis, allowing front-dynamics descriptions beyond bifurcations from stationary fronts. Numerical simulations validate the leading-order delay law and illustrate regimes including heterogeneous speed modulation, attraction and repulsion by stationary fronts for $\alpha<0$, and exploratory dynamics outside the strict assumptions.

\medskip

Throughout this paper, we present the main results from \cite{LvV_thesis} as `Results' rather than `Theorems', as this paper does not contain the proof of those results. All Results presented here correspond to Theorems in \cite{LvV_thesis}, whose proof can be found therein.
}

\subsection{Numerical illustration of travelling front solutions}
\begin{figure}[!h]
    \centering
 \scalebox{0.5}{\includegraphics{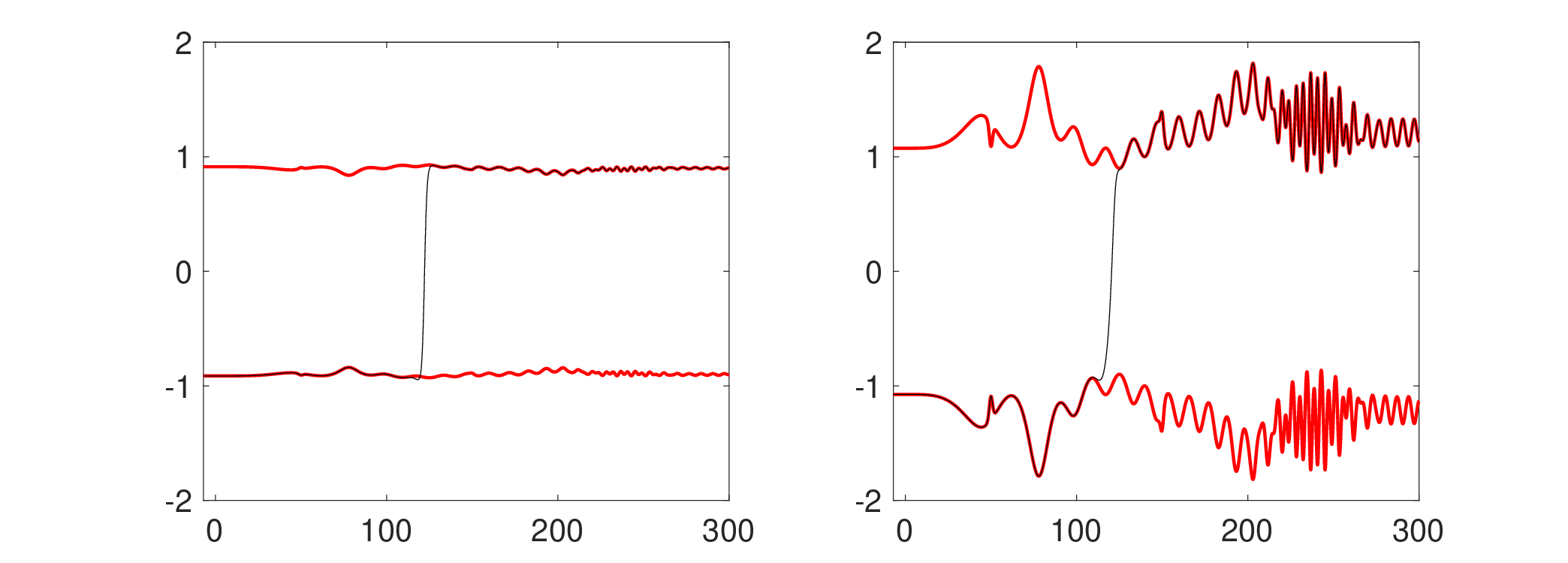}}
\begin{picture}(-100,-100)
    \put(-140,0){$x$}
    \put(80,0){$x$}
    \put(-250,90){\rotatebox{90}{$U(x,t_0)$}}
    \put(-30,90){\rotatebox{90}{$V(x,t_0)$}}
  \end{picture}
\\[1.2cm]
 \scalebox{0.4}{\includegraphics{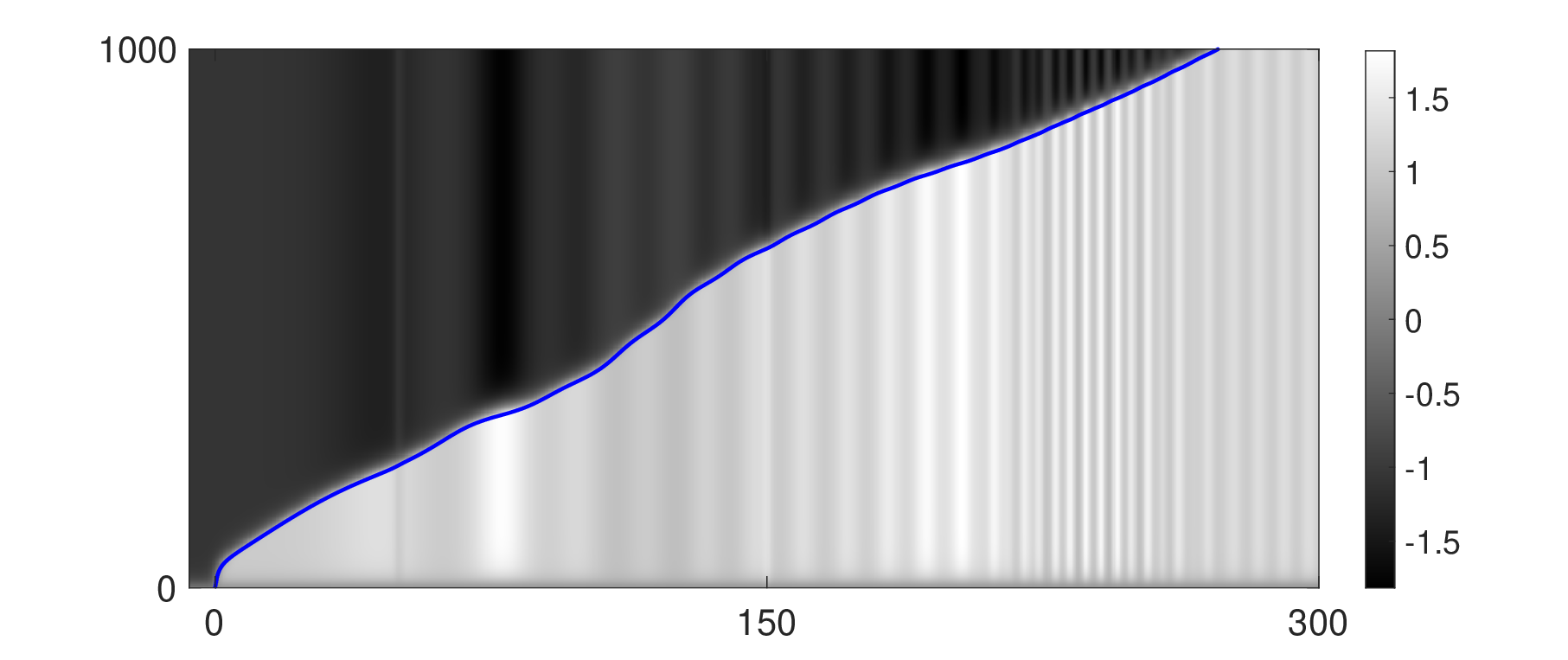}}\\[.5cm]
 \begin{picture}(0,0)
    \put(0,10){$x$}
    \put(-150,105){\rotatebox{90}{$t$}}
    \put(-15,180){$V(x,t)$}
  \end{picture}
 \caption{  Numerical computation of solutions of \eqref{eq:PDE_model}. {\bf Upper panels.} Snapshot in time of the (black) $U$-profile along with (red) background states (left upper panel) and (black) $V$-profile along with (red) background states (right upper panel). {\bf Lower panel.} Density plot for the $V$-profile along with (black) position $z(t)$, defined by the unique zero of the $U$-component, see Definition~\ref{def:position}. {\bf Details of numerical computation.} Coefficients $f_1, f_2$ from \eqref{eq:coefficients_numerics}; Initial Conditions: $u(x,0) = \mathrm{tanh}(x/\varepsilon), v(x,0)=0.1$; Boundary conditions: homogeneous Neumann; Parameter settings: $\varepsilon = 0.15, \alpha = 0.94, \gamma = 0, \hat{\tau} = 1 $. {\it Numerical solver: MATLAB (2020): "pdepe".} Code is available at \cite{chirilusbruckner2026front}}\label{fig:front_profile_numerics}
\end{figure}

Before embarking on a journey to study \eqref{eq:PDE_model} analytically, we briefly show the $U$- and $V-$profiles of travelling front solutions along with a density plot showing the corresponding position $z = z(t)${ , which is defined as the unique zero of the $U$-component (see Definition~\ref{def:position})}. Much of the analytical set-up is, in fact, inspired by careful observations of numerical simulations. In particular, the motionless background states of travelling fronts has been first discovered numerically. In Figure~\ref{fig:front_profile_numerics} we show a traveling front for
\begin{align}\label{eq:coefficients_numerics}
\begin{array}{rcl}
     f_1(x)&=&4+e^{-(x-150)^2}, \qquad f_2(x) \ \ = \ \ 3+0.8g(x)  
    \, ,\\[.2cm]
     g(x)&=&
e^{-(x-50)^2}
-e^{-(0.1(x-80))^2}
+e^{-(0.05(x-120))^2}
-e^{-(0.05(x-200))^2}
+2\,e^{-(0.05(x-240))^2}\cos(1.5x)
+0.5\cos \bigl((0.04x)^2\bigr)
\end{array}
\end{align}

\subsection{Related work and novelty of the present work}
The transition from \eqref{eq:PDE_model_const} to \eqref{eq:PDE_model} is motivated by prior work on a related reaction-diffusion equation, namely, the extended Klausmeier model with spatially varying coefficients, which is used to model the dynamics of vegetation patterns in semi-arid regions where the spatially varying coefficients account for the change in topography (cf. \cite{BC-BD.2020}). This further generalised previous studies which assumed the terrain to be either flat or with a constant slope (\cite{Eric2015}). The novelty in \cite{BC-BD.2020} consisted in using a blend of Geometric Singular Perturbation Theory (GSPT) and the theory of exponential dichotomies to construct spatially varying background states and stationary pulse solutions that are bi-asymptotic to those for $x \rightarrow \pm \infty$. These new types of pulse solutions with heterogeneous tails model vegetation patches on varying topography. Going beyond the stationary case, a reduced description of the position of travelling vegetation patterns was given consisting of a coupled ODE-algebraic system that described the drift instability of stationary pulses. While similar in methodology to \cite{BC-BD.2020}, the study of front solutions of \eqref{eq:PDE_model} in the present work poses several new challenges whose resolution will allow to generalize the approach to a larger class of equations. We will comment on those throughout the text.

\medskip

Going beyond the dynamics of single vegetation patches, the influence of heterogeneity on $N$-pulse (also sometimes called $N$-spike) solutions has been studied by a reduction procedure to a finite-dimensional systems of ODEs, e.g. in \cite{Robbin2019_N_pulse} and \cite{Theo2018}. Furthermore, the limit of small-amplitude vegetation patterns on spatially periodic terrains has been studied in paradigm models such as the Swift-Hohenberg equation \cite{Ehud2008,Jolien2025} by deriving a Ginzburg-Landau equation for the pattern amplitude. In \cite{Ehud2015} even the idea of reversing desertification via resonance with a spatially periodic forcing coefficient was presented.

\medskip

For scalar Nagumo/Allen-Cahn-type settings, in \cite{MR.2024} and \cite{BDK.2025} the influence of heterogeneity on front and multi-front solutions has been studied. Related results can be also found in \cite{AT.2017} for the Allen-Cahn with forcing and discrete versions thereof. Going back to reaction-diffusion systems, in \cite{Peter2012_het} a three-component model of the form \eqref{eq:PDE_model_const} was studied with $\gamma = \gamma(x)$ a step function. All of these settings can be viewed as complementary results to the present work, since they explain how heterogeneity influences primarily the fast variable $U$.  

\medskip

The influence of heterogeneity on pattern formation in 2-D is a very active research field with many applications including to universal phenomena like the formation of Turing patterns or spiral waves. We refrain from giving an extensive list, but rather point to two bodies of work that contain a wealth of analytical results and techniques: \cite{KollarScheel2007}, \cite{JS.2015}, \cite{Ja.2015}, \cite{JSW.2019}, \cite{J.2023} and also \cite{VanGorder2021PatternFormation},\cite{Krause2020WKBJ}, \cite{Krause2025_SH_het}. All these settings have many similarities with \eqref{eq:PDE_model} and will serve as an inspiration for future work.

\newpage

\section{Stationary background states}\label{s:background_states}
To study the existence of front solutions in \eqref{eq:PDE_model}, we first have to determine the existence of stationary background states that can be connected with via a front -- in other words, to which the front is bi-asymptotic as $x \to \pm \infty$. 

\begin{res}\label{res:backgroundstates}
\textcolor{black}{Let $f_1$ and} $f_2$ satisfy Assumptions~\ref{ass:f}. Then there exists $\eps_0 > 0$ such that for all $0 \leq \eps \leq \eps_0$ system \eqref{eq:PDE_model} has a pair of bounded stationary solutions $(u_b^\pm(x;\eps),v_b^\pm(x;\eps))$ that satisfy
\begin{align}
  \inf_{\eps\in [0,\eps_0] }\inf_{x\in \mathbb{R}} u_b^+(x;\eps)> 0  \text{ and } \sup_{\eps\in [0,\eps_0]}\sup_{x\in \mathbb{R}} u_b^-(x;\eps)< 0,  
\end{align}
\begin{align}
  \inf_{\eps\in [0,\eps_0] }\inf_{x\in \mathbb{R}} v_b^+(x;\eps)> 0  \text{ and }   \sup_{\eps\in [0,\eps_0] }\sup_{x\in \mathbb{R}} v_b^-(x;\eps)< 0.  
\end{align}
Furthermore, $u_b^\pm(\cdot;\eps)\in C_b^7(\mathbb{R})$, $v_b^\pm(\cdot;\eps)\in C_b^5(\mathbb{R})$, $u_b^\pm (x;0)= \pm 1$. Suppose $f_1=0$. Using the solution operator 
\begin{equation}\label{eq:solutionoperator_Veq}
    G(\phi)(x) := e^{-x}\int_{-\infty}^x \frac{1}{2}e^\xi \phi(\xi)\,\mathrm{d}\xi + e^x \int_x^\infty \frac{1}{2}e^{-\xi}\phi(\xi)\,\mathrm{d}\xi \, ,
\end{equation}
we have that
\begin{equation}\label{eq:v_b_formula}
 v_b^\pm (x;0) = \pm G(1+f_2)(x) =  \pm e^{-x}\int_{-\infty}^x \frac{1}{2}e^\xi (1+f_2(\xi))\mathrm{d}\xi  \pm  e^x \int_x^\infty \frac{1}{2}e^{-\xi}(1+f_2(\xi))\,\mathrm{d}\xi \, .
\end{equation}
\end{res}

\begin{rmk}
 For a general version of \eqref{eq:v_b_formula} for nonzero $f_1$, see \cite[Chapter 3]{LvV_thesis} \end{rmk}
\begin{rmk}
We observe that if $f_2 = 0$ then $v_b^{\pm} = \pm1$. Furthermore, the result holds equally for non-zero $f_1$ with the small change that \eqref{eq:v_b_formula} will feature a fundamental system of $v'' - (1+f_1(x))v = 0$ instead of exponentials. Notice also that, to leading order, the effect of $f_k$ is only visible in the $v$-component. 
\end{rmk}

Result \ref{res:backgroundstates} is stated as a Theorem in \cite[Chapter 3]{LvV_thesis}, where its full proof can be found. We omit the details of this proof here, and provide an outline of the proof structure. 
The proof of Result \ref{res:backgroundstates} is based on the study of  the ODE system 
\begin{subequations}\label{eq:background_ODEsystem}
\begin{align}
    \eps u_x &= p,\\
    \eps p_x &= -u+u^3 + \eps(\alpha v+\gamma),\\
    v_x &= q,\\
    q_x &= (1+f_1(\chi))v - (1+f_2(\chi))u,\\
    \chi_x &= 1.
\end{align}
\end{subequations}

\noindent We are interested in finding background orbits $\textnormal{BO}_\varepsilon^\pm$ of \eqref{eq:background_ODEsystem} for $0\le \varepsilon \ll 1$ (the singular limit $\varepsilon=0$ included where \eqref{eq:background_ODEsystem} becomes a DAE, i.e. a differential algebraic equation ), which are characterized by  imposing that the $(u,p,v,q)$ components are bounded, and that the $(u,p)$ components are $\mathcal{O}(\epsilon)$ close to $(\pm 1, 0)$.
The usual path for analysis in case $f_k = 0$ would be to exploit the slow-fast structure that gives rise to a critical manifold of which only a compact subset is considered. A similar approach for non-zero $f_k$ leading to the system \eqref{eq:background_ODEsystem} already poses the first difficulty, since one now cannot restrict to compact subsets of the critical manifold, so standard Fenichel theory \cite{Fenichel.1974,Fenichel.1977,Fenichel.1979} cannot be applied directly. Instead, the proof in \cite{LvV_thesis} combines the approach of \cite{Fenichel.1974,Fenichel.1979} with the results on exponential dichotomies found in \cite{Fenichel.1996}. In the singular limit $\eps = 0$, the existence of a pair of (reduced)  background orbits $\textnormal{BO}_{0,R}^\pm$ in the slow reduced system
\begin{subequations}\label{eq:background_ODEsystem_LvV_slowred}
\begin{align}
    v_x &= q,\\
    q_x &= [1+f_1(\chi)]v \pm [1+f_2(\chi)],\\
    \chi_x &= 1
\end{align}
\end{subequations}
for $u = \pm 1$ is shown using the existence of an exponential dichotomy in the homogeneous equation $v_{xx} = [1+f_1]v$, which is equivalent to the Riccati equation $a^2 + a_x = 1 + f_1(x)$. The background orbit for $\varepsilon=0$ is given by $\textnormal{BO}^\pm_{0}:=\{u=\pm 1, p=0\}\times\textnormal{BO}_{0,R}^\pm$.

\begin{rmk}
Note that if the slow variable $v$ were to satisfy a more general linear system one could still use the theory of exponential dichotomies as in \cite{Coppel.1978, BC-BD.2020} with the drawback that one then needs to enforce a (additional) bound on $f_1$ which the Riccati equation for the special case of $v_{xx} = (1+f_1)v$ circumvents. 
The existence of orbits $\text{BO}_0^\pm$ to \eqref{eq:background_ODEsystem_LvV_slowred} then follows, cf. \cite[Prop.2, Ch.8]{Coppel.1978}. 
\end{rmk}

\begin{figure}[!t]
    \begin{overpic}[width=0.9\textwidth]{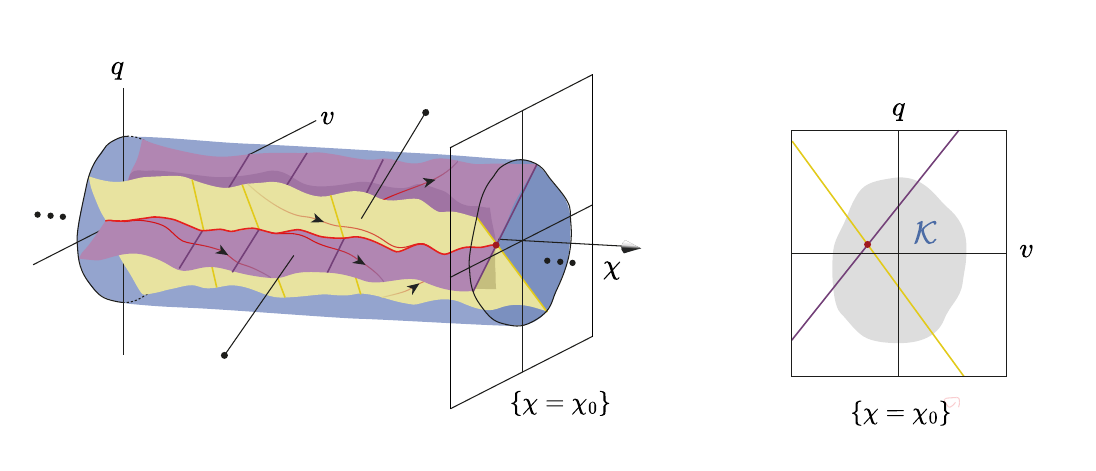}
    \put(50,28){\textcolor{myblack}{$\mathcal{M}_0^+$}}
    \put(12,22){\textcolor{red}{$\text{BO}^+_0$}}
    \put(33,36){\textcolor{myyellow}{$W^{u}_\text{slow}(\text{BO}^+_0)$}}
    \put(12,6){\textcolor{mypurple}{$W^{s}_\text{slow}(\text{BO}^+_0)$}}
    \put(87,25){\textcolor{mypurple}{$\ell^{s,+}(\chi_0)$}}
    \put(87,12){\textcolor{myyellow}{$\ell^{u,+}(\chi_0)$}}
    \end{overpic}
    \caption{A visualisation of the reduced slow dynamics in a neighbourhood of $\text{BO}^+_0$ in the hyperplane $\left\{ u=+1,\, p=0\right\}$. Left: the manifold $\mathcal{M}_0^+$ in black, the singular bounded orbit $\text{BO}^+_0$ in red, the slow unstable manifold of $\text{BO}^+_0$ in purple, the slow stable manifold of $\text{BO}^+_0$ in yellow. Right: the intersection with the hyperplane $\left\{\chi = \chi_0\right\}$, with the line $\ell^{u,+}(\chi_0)$ in purple and $\ell^{s,+}(\chi_0)$ in yellow; the compact set $\mathcal{K}$ is indicated in gray.}\label{fig:T0}
\end{figure}

\noindent Note that, for $f_1 \equiv 0$, the $v$-component of the orbit $\textnormal{BO}_0^\pm$ yields the solution $v_b^\pm(x;0)$ as given in Result \ref{res:backgroundstates}. The persistence of the orbits $\textnormal{BO}_0^\pm$ for $0 \le \eps \ll 1$ in the full system \eqref{eq:background_ODEsystem} is  proved using the following steps. First, we define the manifolds 
\begin{equation}\label{eq:T0pm}
\mathcal{M}_0^\pm := \left\{u=\pm 1,\;p=0, (v,q) \in \mathcal{K},\chi\in \mathbb{R}\right\},
\end{equation}
with $\mathcal{K} \subset \mathbb{R}^2$ a compact set in $(v,q)$-space chosen sufficiently large to (at least) contain the appropriate singular background orbit 
\begin{equation}\label{eq:BO0pm}
\text{BO}^\pm_0 := \{u=\pm 1,p=0\}\times \bigcup_{\chi\in \mathbb{R}}(v_b^\pm(\chi;0), q_b^\pm(\chi;0),\chi) \,, \qquad q_b^\pm = \frac{\text{d} v_b^\pm}{\text{d}x} \, .   
\end{equation} 
The slow stable and unstable manifolds of $\text{BO}^\pm_0$  (defined as the stable  and unstable manifolds of  $\textnormal{BO}_{0,R}^\pm\subseteq \mathcal{M}_0^\pm$) \eqref{eq:background_ODEsystem_LvV_slowred} are foliated by a family of lines $\ell^{u/s,\pm}(\chi)$, see Figure~\ref{fig:T0}, i.e.
\begin{equation}\label{eq:lines}
    W_\text{slow}^{u/s}\left(\text{BO}^\pm_0\right) = \bigcup_{\chi \in \mathbb{R}} \ell^{u/s,\pm}(\chi) \,.
\end{equation}
The manifolds $\mathcal{M}_0^\pm$ are normally hyperbolic, which is immediately apparent from the fast $(u,p)$-dynamics in \eqref{eq:background_ODEsystem}. This allows us to define the (fast) local stable and unstable manifolds of $\mathcal{M}^\pm_0$,
\begin{equation}\label{eq:Wus_T0pm}
    W^{u/s}(\mathcal{M}_0^\pm) = \left\{\,u\in \mathcal{U}^\pm ,\;p = -\frac{1-u^2}{\sqrt{2}}\text{ (stable) or }p = \frac{1-u^2}{\sqrt{2}}\text{ (unstable)},\;(v,q,\chi) \in \mathcal{M}_0^\pm\,\right\}.
\end{equation}

\noindent Here $\mathcal{U}^\pm$ is a compact interval containing $\pm 1$ and not containing $\mp 1$.\\

\begin{figure}[!t]
    \begin{overpic}[width=0.9\textwidth]{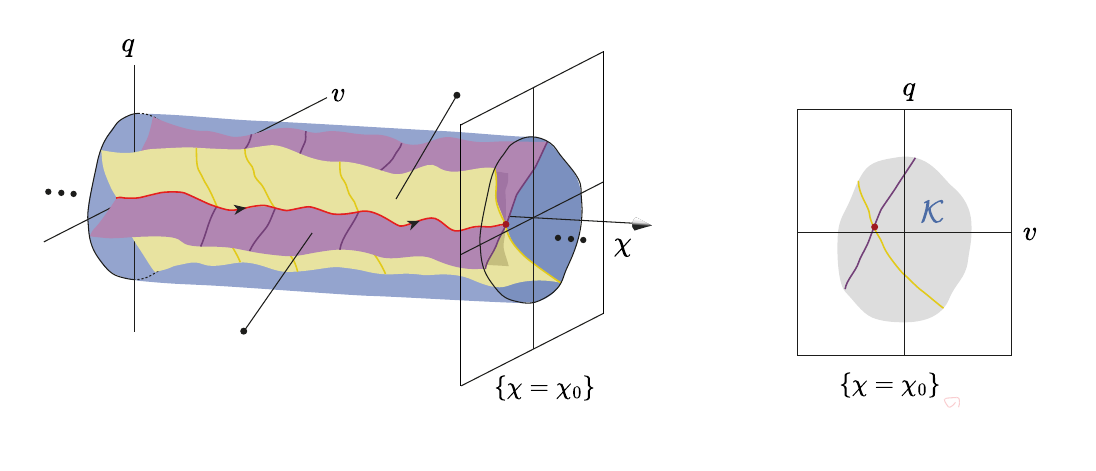}
    \put(50,31){\textcolor{myblack}{$\mathcal{M}_\eps^+$}}
    \put(12,22){\textcolor{red}{$\text{BO}^+_\eps$}}
    \put(35,37){\textcolor{myyellow}{$W^{u}_\text{slow}(\text{BO}^+_\eps)$}}
    \put(13,8){\textcolor{mypurple}{$W^{s}_\text{slow}(\text{BO}^+_\eps)$}}
    \put(85,29){\textcolor{mypurple}{$\mathcal{F}^{s,+}(\chi_0)$}}
    \put(87,13){\textcolor{myyellow}{$\mathcal{F}^{u,+}(\chi_0)$}}
    \end{overpic}
    \caption{A visualisation of the slow dynamics in a neighbourhood of $\text{BO}^+_\eps$ near the hyperplane $\left\{\,u=+1,\;p=0\,\right\}$. Left: the manifold $\mathcal{M}_\eps^+$ in black, the background orbit $\text{BO}^+_\eps$ in red, the slow unstable manifold of $\text{BO}^+_\eps$ in purple, the slow stable manifold of $\text{BO}^+_\eps$ in yellow. Right: the intersection with the hyperplane $\left\{\chi = \chi_0\right\}$, with the fibre $\mathcal{F}^{u,+}(\chi_0)$ in purple and $\mathcal{F}^{s,+}(\chi_0)$ in yellow; the compact set $\mathcal{K}$ is indicated in gray.}\label{fig:Teps}
\end{figure}

The main technical result of \cite[Chapter 3]{LvV_thesis} is that $W^{u/s}(\mathcal{M}_0^\pm)$ $C^3$-smoothly perturb for $0<\eps \ll 1$ to $W_\eps^{u/s,\pm}$. Their transversal intersection defines $\mathcal{M}_\eps^\pm := W_\eps^{u,\pm} \cap W_\eps^{s,\pm}$, which is a $C^3$-smooth perturbation of $\mathcal{M}_0^\pm$. Moreover, $W_\eps^{u/s,\pm}$ are locally invariant unstable and stable manifolds of $\mathcal{M}_\eps^\pm$ in the sense of Fenichel -- that is, orbits in $W_\eps^{u,\pm}$ approach $\mathcal{M}_\eps^\pm$ backwards in $x$ with exponential rate, and orbits in $W_\eps^{s,\pm}$ approach $\mathcal{M}_\eps^\pm$ forwards in $x$ with exponential rate. The proof of this result follows the steps taken in \cite{Fenichel.1979}, in particular using the Hadamard graph transform, where uniform estimates obtained therein from the compactness of the underlying critical manifold are replaced by uniform estimates on the vector field \eqref{eq:background_ODEsystem}.  The proof of smoothness in $\eps$ (see \cite[Chapter 3]{LvV_thesis} for full details) follows the approach of \cite{Fenichel.1974}, in particular Theorem 6 therein. These persistence results provide us with the existence of a function $u^\pm(v,q,\chi;\eps) = \pm 1 + \mathcal{O}(\eps)$ parametrising the $u$-component of $\mathcal{M}^\pm_\eps$ as a graph over $(v,q,\chi)$, which can be used to study the slow dynamics on $\mathcal{M}_\eps^\pm$, given by
\begin{subequations}
\label{eq:background_ODEsystem_LvV_slow_eps}
\begin{align}
    v_x &= q,\\
    q_x &= (1+f_1(\chi))v - (1+f_2(\chi))u^\pm(v,q,\chi;\eps),\\
    \chi_x &= 1.
\end{align}
\end{subequations}
The results on nonlinear exponential dichotomies by Fenichel \cite{Fenichel.1996} can now be used to prove that the slow stable and unstable manifolds $W^{u/s}_{slow}(\text{BO}^\pm_0)$ perturb $C_1$-smoothly (w.r.t. the dynamics \eqref{eq:background_ODEsystem_LvV_slow_eps}), and that their transversal intersection $\textnormal{BO}_{\varepsilon,R}^\pm$ is a $C_1$-smooth perturbation of the reduced singular bounded orbit $\text{BO}^\pm_{0,R}$. 
In addition, the (perturbed) manifolds $W^{u}_\text{slow}(\text{BO}^\pm_\eps)\subseteq \mathcal{M}_\varepsilon^\pm$ and $W^{s}_\text{slow}(\text{BO}^\pm_\eps)\subseteq  \mathcal{M}_\varepsilon^\pm$ are the stable and unstable manifolds of  $\textnormal{BO}_{\varepsilon,R}^\pm$ in the sense of Fenichel (w.r.t. the dynamics \eqref{eq:background_ODEsystem_LvV_slow_eps}). We obtain the background orbits $\textnormal{BO}_{\varepsilon}^\pm$ in \eqref{eq:background_ODEsystem} as the image of $\textnormal{BO}_{\varepsilon,R}^\pm$ under the inclusion map $\mathcal{M}_\varepsilon^\pm\xhookrightarrow{}\mathbb{R}^5$.
The resulting geometry is similar to that in the singular limit; now, the intersection of $W^{u/s}_\text{slow}(\text{BO}^\pm_\eps)$ with the hyperplane $\left\{\chi = \chi_0\right\}$ is a Fenichel fibre $\mathcal{F}^{u/s,\pm}(\chi_0)$, see also Figure \ref{fig:Teps}. The $(u,v)$-components of the perturbed background orbits $\text{BO}^\pm_\eps$ correspond to the bounded stationary background states whose existence and properties are stated in Result \ref{res:backgroundstates}.

\begin{rmk}
Already foreshadowing the construction of stationary front solutions via Melnikov analysis, we want to emphasise that in the compact case it is sufficient to have that $W^{u/s}(\mathcal{M}_0^\pm)$ is a $C^2$-smooth perturbation to ensure $C^1$-smoothness of the Melnikov function. In the non-compact case one needs to go a step further and demand bounded $C^3$-smoothness which then implies uniform $C^2$-smoothness without compactness.
\end{rmk}

\begin{rmk}
In \cite{LvV_thesis} is it shown that the background solutions from Result~\ref{res:backgroundstates} are stable, hence, making \eqref{eq:PDE_model} a bi-stable PDE.
\end{rmk}

\begin{rmk}
In the present equation setting the critical and slow manifold are different. This is in contrast with the corresponding system in \cite{BC-BD.2020} where the critical and slow manifold coincide. In that respect, the setting studied here is more representative of the general case.   
\end{rmk}

\begin{figure}[H]
\centering
    \begin{overpic}[width=0.6\textwidth]{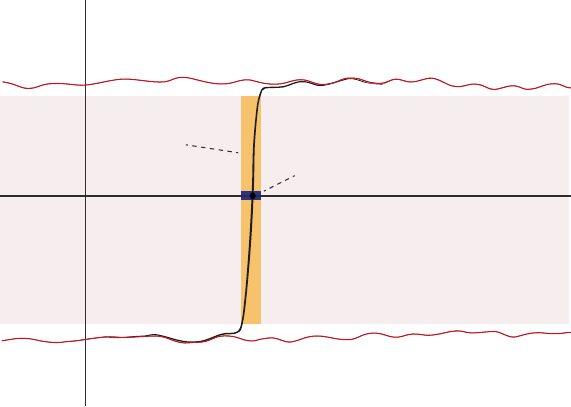}
    \put(50,42){\textcolor{mydarkblack}{$I_t = U_{\text{TF}}(\cdot,t)^{-1}(J)$}}
    \put(11,16){\textcolor{mydarkblack}{$u_l$}}
    \put(11,51){\textcolor{mydarkblack}{$u_r$}}
    \put(28,45){\textcolor{myorange}{$J$}}
    \put(16,70){$u$}
    \put(99,33){$x$}
    \put(102,55){\textcolor{mydarkred}{$u_b^+(x)$}}
    \put(102,12){\textcolor{mydarkred}{$u_b^-(x)$}}
    \end{overpic}
    \caption{A visual representation of a travelling front solution as described in Definition \ref{def:front}.}\label{fig:frontdef}
\end{figure}
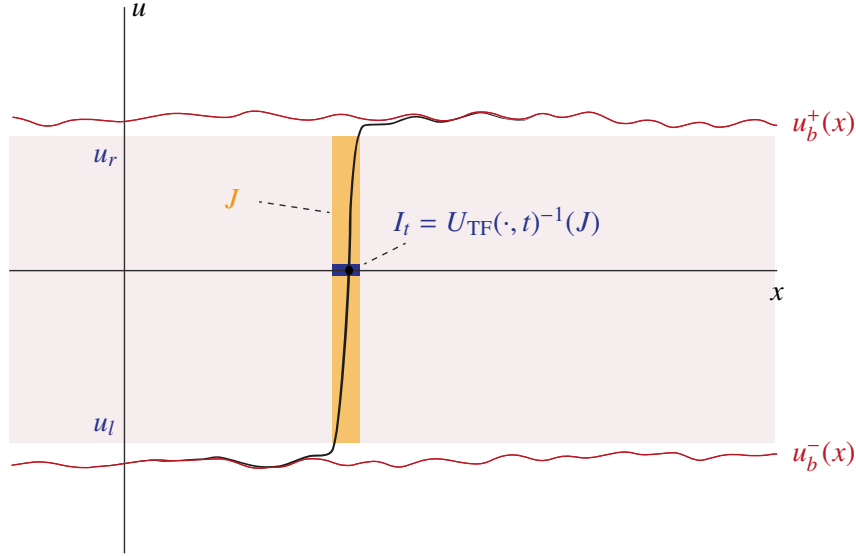

\newpage

\section{Front solutions}\label{s:front_solutions}
While the concept of a travelling front is intuitively unambiguous, it turns out to be useful to provide a formal description of a travelling front solution to \eqref{eq:PDE_model}, as the presence of spatial heterogeneities introduces non-trivial spatio-temporal features to such solutions. 

\begin{deftn}[\bf Front solutions]\label{def:front}
 Let $(u_b^\pm,v_b^\pm)$ be the stationary background states from Result \ref{res:backgroundstates}. A front solution $\left(U_\text{F},V_\text{F}\right)$ is a bounded solution of \eqref{eq:PDE_model} satisfying the following properties for all $t \in \Psi $ where either $ \Psi = [t_0, \infty)$ for some $t_0 \in \mathbb{R}$ or $\Psi = \mathbb{R}$:
 \begin{enumerate}
     \item \textbf{(Regularity)} $\left(U_\text{F},V_\text{F}\right)$ is $C^2$-smooth in $x$ and $C^1$-smooth in $t$ for $t \in \Psi$.
     \item \textbf{(Asymptotic behaviour in space)} For every $t \in \Psi$, $\left(U_\text{F}(\cdot,t),V_\text{F}(\cdot,t)\right)$ is asymptotic to the stationary background states $(u_b^\pm(x),v_b^\pm(x))$ as $x \to \pm \infty$.
     \item \textbf{(Front interface properties)} There exists an interval $J := [u_l,u_r]$, $u_l<0<u_r$ such that for all $t \in \Psi$
     \begin{enumerate}
         \item $I_t := U_{\text{F}}^{-1}(\cdot,t)\left(J\right)$ is an interval, \label{def:front_fastint}
         \item $U_{\text{F}}(\cdot,t)$ is a bijection between $I_t$ and $J$, \label{def:front_monotone}
         \item $0 < \text{inf}_{t \in \Psi} \left| I_t \right|$ and $\text{sup}_{t \in \Psi} \left| I_t \right| < \infty$\textcolor{black}{, where $\left| I_t\right|$ is the length of the interval $I_t$}. \label{def:front_time_uniform}
     \end{enumerate}
     \item \textbf{(Uniform decay in tails)} The decay of the solution at $\pm \infty$ to the background states is uniform in $t \in \Psi$.
 \end{enumerate}
\end{deftn}

\begin{deftn}[\bf Entire front solutions]\label{def:entire_front}
 We call a front solution from Definition~\ref{def:front} an entire front if $\Psi = \mathbb{R}$.
\end{deftn}

The front interface properties described in Definition \ref{def:front} characterise a travelling front solution by the fact that a clearly defined front connection between the two stationary background states persists for an open time interval $\Psi \subset \mathbb{R}$. In particular, properties \ref{def:front_fastint} and \ref{def:front_monotone} specify that, for fixed $t$, this front \textcolor{black}{interface} is a monotonically increasing graph over an $x$-interval $I_t$. Property \ref{def:front_time_uniform} ensures that this interval persists in time, i.e. does not disappear or stretch to infinity. For a visual representation of Definition \ref{def:front}, see Figure \ref{fig:frontdef}.Note that Definition \ref{def:front} confines our analysis to fronts that are left asymptotic to $(u_b^-,v_b^-)$ and right asymptotic to $(u_b^+,v_b^+)$, for which the connecting front is monotonically increasing in $x$. This choice has been made for clarity of presentation; fronts that are left asymptotic to $(u_b^+,v_b^+)$ and right asymptotic to $(u_b^-,v_b^-)$ can be analysed in an analogous manner.

\medskip
While in the constant-coefficient case \eqref{eq:PDE_model_const} it is taken for granted that uniformly travelling front solutions exist for all time, a similar result for \eqref{eq:PDE_model} will require considerably more effort. In particular, it will turn out to be necessary to first treat the initial value problem (leading to Definition~\ref{def:Initialised_front} where $\Psi$ from Definition~\ref{def:front} is not $\mathbb{R}$) to then prove the existence of entire solutions as given in Definition~\ref{def:entire_front}.

\subsection{Stationary front solutions}\label{ss:stationaryfronts}
Definition \ref{def:front} does not specify any temporal behaviour of travelling front solutions to \eqref{eq:PDE_model}: in particular, stationary front solutions fall within the class of solutions described in Definition \ref{def:front}. Although the main focus of this paper is on the existence of genuinely travelling fronts, i.e. those with non-trivial temporal behaviour, it is necessary to briefly treat the existence of stationary front solutions to \eqref{eq:PDE_model}, as these can serve as attractors for the front solutions to \eqref{eq:PDE_model} that do evolve in time.

\begin{res}[\bf 1-parameter family of stationary fronts]\label{res:stationaryfronts} 
Let the conditions of Result~\ref{res:backgroundstates} be fulfilled such that there exist background states $(u_b^{\pm}, v_b^{\pm})$. There exists $\widetilde{\varepsilon}_0 > 0$ ($\varepsilon_0$ from Result~\ref{res:backgroundstates}) and a unique $C^1$-smooth function $\gamma_{SF}(x_0;\varepsilon)$ defined
on $\mathbb{R}\times [0,\varepsilon_0]$ such that for $\varepsilon\in (0,\varepsilon_0]$ and $\gamma = \gamma_{SF}(x_0;\varepsilon)$ there exists a stationary front as defined in Definition~\ref{def:front} of the PDE~(\ref{eq:PDE_model}). Its spatial profile
	$(U_{SF},V_{SF})$ (which depends on the parameters $x_0,\varepsilon$)  satisfies the following properties:
	\begin{enumerate}
		\item \textbf{(Front position)} $U_{SF}(x_0;x_0, \varepsilon)=0$.
	
		\item \textbf{(Regularity and limiting behaviour)} The $v$-component of the front  $(x,x_0,\varepsilon)\to V_{SF}(x;x_0, \varepsilon)$ extends $C^1$-smoothly to $\{\varepsilon=0\}$. The $u$-component of the front $(x,x_0,\varepsilon)\to U_{SF}(x;x_0, \varepsilon)$ extends $C^1$-smoothly to $\{\varepsilon=0,x\neq x_0 \}$ and its limit satisfies $U_{SF}(x;x_0, 0)= -1+2 \chi_{ ( x_0,\infty]  }(x)$.
	\item  \textbf{(Leading term of existence condition)}  
		$\gamma_{SF}(x_0;0)$ is given by 
        \begin{eqnarray}
		\alpha v_{SF}(x_0;0)+ \gamma_{SF}(x_0;0)=0,\label{gammaSF0}   
		\end{eqnarray}
		$v_{SF}(x_0;0):=V_{SF}(x_0;x_0,0)$ is given by
		\begin{eqnarray}
		v_{SF}(x_0;0)= \frac{2\frac{d}{dx}v_b^-(x_0;0)-v_b^-(x_0;0)(a_s(x_0)+a_u(x_0)) }{ a_s(x_0;
			)- a_u(x_0)} ,\label{vSF}
		\end{eqnarray}
		\noindent arising from the intersection of lines $\ell^{u,-}(x_0)$ and $\ell^{s,+}(x_0)$ from \eqref{eq:lines} that have slopes 
		$a_{s/u}(x)$ given by unique positive, respectively negative bounded solutions of the Ricatti differential equation
			\begin{eqnarray} 
			a_x= 1+ \delta_1\widetilde{f}_1(x)-a^2\, , \quad f_1 = \delta_1 \widetilde{f}_1 \, .
			\end{eqnarray}  
		Furthermore, one has $\|a_{u}(.)- 1\|_\infty=O(\delta_1)$ and $\|a_s(.)+1\|_\infty = O(\delta_1)$.	
	\end{enumerate} 
\end{res} 

Result \eqref{res:stationaryfronts} is stated as a Theorem in \cite[Chapter 3]{LvV_thesis}, where its full proof can be found. As with Result \ref{res:backgroundstates}, we omit the details of that proof here, and provide a brief outline of the proof structure. \\
The proof of Result \ref{res:stationaryfronts} uses the geometric concepts that were introduced in the proof of Result \ref{res:backgroundstates}, in particular the stable and unstable manifolds $W^{u/s,\pm}_\eps$ of the manifold $\mathcal{M}^\pm_\eps$, which to leading order in $\eps$ are given by $W^{u/s}(\mathcal{M}_0^\pm)$ \eqref{eq:Wus_T0pm} and $\mathcal{M}^\pm_0$ \eqref{eq:T0pm}, respectively. A stationary front solution to \eqref{eq:PDE_model} is a heteroclinic orbit in \eqref{eq:background_ODEsystem} that connects $\text{BO}_\eps^-$ to $\text{BO}_\eps^+$; these background orbits are given to leading order in $\eps$ by $\text{BO}_0^\pm$ \eqref{eq:BO0pm}.\\
A heteroclinic orbit from $\text{BO}_\eps^-$ to $\text{BO}_\eps^+$ lies in the intersection $W^{u,-}_\eps \cap W^{s,+}_\eps$, which to leading order in $\eps$ is given by 
\begin{equation}\label{eq:Wu-_int_Ws+}
W^u(\mathcal{M}_0^-) \cap W^s(\mathcal{M}_0^+) = \left\{\,u\in \mathcal{U}^+\cap \mathcal{U}^- ,\;p = \frac{1-u^2}{\sqrt{2}},\;(v,q,\chi) \in \mathcal{M}_0^\pm\,\right\}.
\end{equation}
In the fast spatial coordinate $\xi = \frac{x}{\eps}$, the relation $u_\xi = p = \frac{1}{\sqrt{2}}(1-u^2)$ \eqref{eq:background_ODEsystem} yields the leading order front profile $\text{tanh} \frac{\xi}{\sqrt{2}}$. The transversality and persistence in $\eps$ of the intersection \eqref{eq:Wu-_int_Ws+} is proven in \cite{LvV_thesis} by the study of the properties of a Melnikov function $M(v,q,\chi;\eps)$ that measures the distance between $W^{u,-}_\eps$ and $W^{s,+}_\eps$, which to leading order in $\eps$ is given by 
\begin{equation}\label{eq:stat_Melnikov}
M(v,q,\chi;0) = 2\sqrt{2}(\alpha v + \gamma). 
\end{equation}
Furthermore, the dynamics on $\mathcal{M}^\pm_\eps$, which are given by \eqref{eq:background_ODEsystem_LvV_slow_eps}, are characterised by $W^{u/s}_\text{slow}(\text{BO}_\eps^\pm)$, see also Figure \ref{fig:Teps}. To leading order in $\eps$, we have
\begin{equation}
    W^{u/s}_\text{slow}(\text{BO}_0^\pm) = \left\{u = \pm 1,\,p=0,\,q = \frac{\text{d}v_b^\pm}{\text{d}x}(\chi;0) + a_{u/s}(\chi)\left(v-v_b^\pm(\chi;0)\right)\right\},
\end{equation}
where $a_{u/s}$ are the slopes of the unstable/stable lines corresponding to the exponential dichotomy of the homogeneous equation $v_{xx} = (1+f_1(\chi))v$. As the stationary front described in Result \ref{res:stationaryfronts} lies in the intersection $W^u(\text{BO}_\eps^-) \cap W^s(\text{BO}_\eps^+)$, we can track their intersection through $W^u_\text{slow}(\text{BO}_\eps^-)$ via the fast transition $W^{u,-}_\eps \cap W^{s,+}_\eps$ to $W^s_\text{slow}(\text{BO}_\eps^+)$. This yields, to leading order in $\eps$, the existence condition
\begin{equation}
\frac{\text{d}v_b^-}{\text{d}x}(\chi;0) + a_u(\chi)\left(v-v_b^-(\chi;0)\right) = \frac{\text{d}v_b^+}{\text{d}x}(\chi;0) + a_s(\chi)\left(v-v_b^+(\chi;0)\right),
\end{equation}
which can be solved for $v$ to give
\begin{equation}
    v_{\rm SF} = \frac{1}{a_s(\chi) - a_u(\chi)} \left(\frac{\text{d}v_b^-}{\text{d}x}(\chi;0) - \frac{\text{d}v_b^+}{\text{d}x}(\chi;0) + a_s(\chi) v_b^+(\chi;0) - a_u(\chi) v_b^-(\chi;0)\right) \, ,
\end{equation}
which together with the leading order Melnikov condition $M(v,q,\chi;0) = 0$ \eqref{eq:stat_Melnikov} yields existence condition \eqref{gammaSF0}.

\begin{rmk} For $f_1 = f_2 = 0$ the leading order existence condition \eqref{gammaSF0} reduces to \eqref{eq:existence_const_coeff}. Moreover, if $f_1 = 0$ and $f_2$ is non-zero fulfilling Assumption~\ref{ass:f} then the existence \eqref{gammaSF0} is given by $
\alpha \frac{d}{dx} v_b^-(x_0;0)+ \gamma_{SF}(x_0;0)=0$. Note that these conditions reflect, in particular, the loss of translation-invariance in space, that is, while in the constant coefficient case the fronts come in families parametrised by an arbitrary position, in the heterogeneous case only specific positions $x_0$ are selected.
\end{rmk}


\subsection{Travelling fronts with non-zero speed}\label{ss:travelingfronts}
We now move to travelling front solutions with with non-zero speed. We will formulate all results for the special case $f_1 = 0$. Similarly as for stationary fronts, the case $f_1 \neq 0$ is only more technical but does not need new techniques. A first step in the construction of travelling entire fronts is a solid command of the initial value problem for travelling fronts. Since this is of interest in its own right, we give a separate definition of {\it initialised fronts}.


\begin{deftn}[\bf Initialised front solutions]\label{def:Initialised_front}
We call a front solution $(U_F, V_F)$ from Definition~\ref{def:front} with initial condition $(U_F, V_F)(x, t_0) = (U_0, V_0)(x)$ for some $t_0 \in \mathbb{R}$ an initialised front if $\Psi = [t_0, \infty)$.
\end{deftn}

\subsubsection{Initialised travelling fronts}\label{ss:TF_IVP}
We now set out to construct front solutions as in Definition~\ref{def:front} with a non-uniformly moving interface connecting {\it motionless} background states. We will work with the fast-reaction scaling
\eqref{eq:PDE_model_fast_reaction} which is equivalent to \eqref{eq:PDE_model}.
The key observation its that, if such solutions exist then for every fixed $s \geq s_0$ their time derivatives \textcolor{black}{satisfy} $(U_s, V_s)(x,s) \rightarrow 0, x \rightarrow \pm \infty$. So, instead of working directly with the solutions $(U, V)$ we instead will study the time evolution of the time derivatives $(U_s, V_s)$. For sufficiently smooth functions $\left(U_0,V_0\right) = \left(U_0,V_0\right)(x)$, we define \textcolor{black}{their generalised `time derivative'} operators through \textcolor{black}{the right-hand side of the vector field of \eqref{eq:PDE_model_fast_reaction} as follows:} 
\begin{subequations}
    \begin{align}
        U_0^{(1)}(U_0, V_0) &:= \frac{1}{\eps^2}\left[\eps^2 \partial_x^2 U_0 + U_0 - U_0^3 -\eps (\alpha V_0 +\gamma)\right],\\
        V_0^{(1)}(U_0, V_0)  &:= \frac{1}{\hat{\tau}}\left[\partial_x^2 V_0 - V_0 + (1+f_2(x))U_0\right].
    \end{align}
\end{subequations}
Note that, if $(U, V) = (U, V)(x,s)$ is a solution of \eqref{eq:PDE_model_fast_reaction} with $(U, V)(x,s_0) = \left(U_0,V_0\right)(x)$ then trivially $(\partial_s U, \partial_s V)(x, s_0) = \left(U_0^{(1)}(U_0, V_0), V_0^{(1)}(U_0, V_0)\right)(x,s_0)$. We will often suppress the arguments $(U_0, V_0)$ in the following for brevity and simply write $U_0^{(1)}$ and $V_0^{(1)}$.

\paragraph{Fronts sets $\Gamma$}
Next, we introduce the vector of positive bounds $\theta = (\theta_1, \theta_2) $ with $ \theta_1:= \left(\overline{u}, \overline{v},u_\ell,u_r\right)$, $ \theta_2:= \left(m_+,m_-,n_+,n_-\right)$. Here $\overline{u}$ and $\overline{v}$ are a-priori-bounds on the solution (details in \cite[Chapter 4]{LvV_thesis}) and $-1<u_\ell<0<u_r<1$ . We
introduce the set $\Gamma^{\rightarrow}(\theta)$ consisting of functions $(U_0,V_0) = (U_0,V_0)(x)$ with the following properties:
\begin{enumerate}
    \item $U_0 \in C^{4+\ell}(\mathbb{R}), V_0 \in C^{4+\ell}(\mathbb{R})$, $\ell\in(0,1)$, \label{prop:regularity}
    \item $|U_0(\cdot)| \leq \overline{u}$, $|V_0(\cdot)| \leq \overline{v}$, \label{prop:apriori_bounds}
    \item $\lim_{x \to \pm \infty} \left(U_0^{(1)}(U_0, V_0)(x), V_0^{(1)}(U_0, V_0)(x)\right) = (0,0)$, \label{prop:asymptotic_timederivative}
    \item The interval $J := (u_\ell,u_r)$ contains $\left[-\frac{1}{\sqrt{3}}, \frac{1}{\sqrt{3}}\right]$ and $I := U_0^{-1}(J)$ is a finite interval, \label{prop:front_existence}
    \item $\partial_x U_0(x) > 0$ for all $x \in I$ ,\label{prop:front_monotonicity}
    \item $\lim_{x \to \pm \infty} |U_0(x) - u_b^\pm(x)| = 0$ and $\lim_{x \to \pm \infty} |V_0(x) - v_b^\pm(x)| = 0$,\label{prop:asymptotic_backgroundstates}
    \item $\left|\mathrm{min}_{x \in \mathbb{R}}\left\{U_0^{(1)}(U_0, V_0)(x),0\right\}\right| \leq \frac{1}{\eps} m_-$ and $ \left|\mathrm{max}_{x \in \mathbb{R}}\left\{U_0^{(1)}(U_0, V_0)(x),0\right\}\right| \leq \eps m_+$, \label{prop:epsbounds_U}
    \item $\left|\mathrm{min}_{x \in \mathbb{R}}\left\{V_0^{(1)}(U_0, V_0)(x),0\right\}\right| \leq n_-$ and $ \left|\mathrm{max}_{x \in \mathbb{R}}\left\{V_0^{(1)}(U_0, V_0)(x),0\right\}\right| \leq \eps n_+$, \label{prop:epsbounds_V}
    \item   $U^{(1)}_0(U_0,V_0)(x)<0$ for $x\in I$.
    \label{directioninterface}
\end{enumerate}

\noindent Likewise we define $\Gamma^\leftarrow(\theta)$, interchanging the role of $m^\pm$ and $n^\pm$ in items 7 (e.g. using $m_-\epsilon$ and $\frac{m^+}{\epsilon}$ in 7.) and 8, and requiring $U^{(1)}_0(U_0,V_0)(x)<0$ for $x\in I$ in 9.
The sets $\Gamma^{\rightarrow}(\theta)$ and $\Gamma^{\leftarrow}(\theta)$ \textcolor{black}{contain} sufficiently smooth bounded function pairs $(U_0(x),V_0(x))$ that are bi-asymptotic to the background states $(u_b^\pm,v_b^\pm)$. The high regularity in 1 is to ensure that the time derivatives $(U_s,V_s)$ of the PDE solution $(U,V)$ are themselves (classical) solutions of the PDE~\ref{eq:PDE_model_UtVt} (shown below). Furthermore, functions in $\Gamma^{\rightarrow}(\theta)$ and $\Gamma^{\leftarrow}(\theta)$ can be characterised as fronts through the front interface properties stated in Definition \ref{def:front}. Their implicitly defined `time derivative' $(U_0^{(1)},V_0^{(1)})$ through the vector field of the system of evolution equations \eqref{eq:PDE_model} converges to zero as $x \to \pm \infty$ and obeys some explicitly $\eps$-dependent bounds that are reminiscent of the scaling relations for conventional uniformly travelling fronts, \textcolor{black}{in contrast to} the scaling relation between derivatives for parabolic equations. \textcolor{black}{Note that the choice of $J$ to contain $\pm \frac{1}{\sqrt{3}}$ (front property 4) is mainly technical, and to some extent arbitrary: the factor $\frac{1}{\sqrt{3}}$ is used in estimates on the size of the difference $|u_r - u_\ell|$.} \\

\noindent One can show that there exist $\theta$ such that the union of the sets $\Gamma^\rightarrow(\theta)$ and $\Gamma^\leftarrow(\theta)$ along with stationary fronts (discussed in Result~\ref{res:stationaryfronts}) is nonempty for $0<\varepsilon \ll 1$. We will briefly discuss the geometric construction for cooking up elements in the proof sketch for Result~\ref{res:TF_noSF} (using auxiliary compactly supported function $\Lambda$). More details can be found in \cite[Chapter 4]{LvV_thesis}.

\paragraph{Temporal evolution of $(\partial_s U, \partial_s V)$.}
As alluded to, it will turn out to be useful to consider the linear PDE system that governs the dynamics of $\partial_s U$ and $\partial_s V$, given by
    \begin{equation}\label{eq:PDE_model_UtVt}
    \left\{
    \begin{aligned}
    \eps^2\partial_s \left(\partial_s U\right) &= \eps^2 \partial_x^2 \left(\partial_s U\right) + (1-3U^2)\partial_s U -\eps \alpha \,\partial_s V, \\
    { \hat{\tau}}\partial_s \left(\partial_s V\right) &= \partial_x^2 \left(\partial_s V\right) - \partial_s V + (1+f_2(x))\partial_s U.
    \end{aligned}
    \right. 
    \end{equation}
System \eqref{eq:PDE_model_UtVt} can be used to determine the signs of $\partial^2_s U$ and $\partial^2_s V$ using similar arguments as in earlier steps of the proof. A crucial part is writing $\partial_s U$ as the sum of three parts: a `fade out' part, which decreases exponentially in time; an `interface' part, whose properties are determined by the shape of the front interface, and a `tails' part, whose properties are determined by the asymptotic convergence to the background states. For all three parts, the variation of constants formula applied to \eqref{eq:PDE_model_UtVt} is used to establish estimates on $|\partial_s U|$.

\paragraph{Positive feedback $\alpha < 0$}
The aforementioned generic regularity results form the basis for 
the next step of the proof. When $\alpha<0$, the set of functions $(U_0,V_0)$ obeying properties \ref{prop:regularity}, \ref{prop:apriori_bounds}, \ref{prop:asymptotic_timederivative}, \ref{prop:asymptotic_backgroundstates} and the additional property
\begin{enumerate}
    \item[10.] $u_b^-(x) \leq U_0(x) \leq u_b^+(x)$ and $v_b^-(x) \leq V_0(x) \leq v_b^+(x)$ \label{prop:bounded_between_backgroundstates}
\end{enumerate}
can be shown to be forward invariant under the flow of \eqref{eq:PDE_model}, using integral estimates combined with repeatedly extending the length of finite time interval on which the results are shown to hold. In the next step, the front properties \ref{prop:front_existence} and \ref{prop:front_monotonicity} are incorporated, together with the propagation direction property
\begin{enumerate}
    \item[11.] $U_0^{(1)}(x,0), V_0^{(1)}(x,0) \leq 0$ for all $x \in \mathbb{R}$ and $U_0^{(1)}(x,0), V_0^{(1)}(x,0) < 0$ for $x \in I=U_0^{-1}([u_\ell,u_r])$. \label{prop:propagation_direction}
\end{enumerate}
It can be shown that the these properties are forward invariant under the flow of \eqref{eq:PDE_model}. 

\paragraph{Mixed feedback $\alpha > 0$}
When $\alpha > 0$, the feedback of $V$ on $U$ in \eqref{eq:PDE_model} is negative, while the feedback of $U$ on $V$ remains positive. This mixed feedback has profound impact on the dynamics of travelling fronts in \eqref{eq:PDE_model}. Indeed, when $\alpha>0$, the propagation direction of the front might reverse; see also Example 2 in section \ref{ss:examples}. A related phenomenon is that, for $\alpha>0$, the $U$-component might overshoot the background states $u_b^\pm$, in contrast to the travelling front described in Results \ref{res:delay_eq} and \ref{res:TF_withSF}.\\

We will now set out to formulate an existence result for $\alpha > 0$ for front solutions that travel in one fixed direction. We formulate two results on the existence of initialised travelling fronts in the absence of stationary fronts. The results are formulated for right-travelling fronts for clarity of presentation; results on left-travelling front can be obtained in an analogous manner.

\begin{res}\label{res:inTF_compactsupport}
Let $\alpha >0$ and $\sup_{x \in \mathbb{R}} \frac{\mathrm{d}v_b^-}{\mathrm{d}x}(x;0) < \frac{\gamma}{\alpha}$. Fix $\theta_1 = (\overline{u},\overline{v},u_\ell,u_r)$, and $\theta_2 = (m_+,m_-,n_+,n_-) \in \mathbb{R}_{>0}^4$. Then there exists $\widehat{\theta}_2$ and $\varepsilon_0 > 0$ such that for all $\varepsilon \in (0, \varepsilon_0)$ and all initial conditions $(U_0, V_0) \in \Gamma^\rightarrow(\theta_1, \theta_2)$ such that $U_0(z_0) = 0$ and $U_0^{(1)}$ \textcolor{black}{fulfils the condition
\begin{align}
   \text{there exist } \varkappa, \lambda_0, \underline{d} > 0 \text{ such that } \text{supp } U^{(1)}(U_0, V_0) \subset [z_0 - \varepsilon \varkappa, z_0 + \varepsilon \varkappa]\\ \text{and } U_0^{(1)} \leq \frac{\lambda_0}{\varepsilon} \text{ for } x \in \left[U_0^{-1}\left(-1/\sqrt{3}\right)- \varepsilon \underline{d}, U_0^{-1}\left(1/\sqrt{3}\right) + \varepsilon \underline{d}\right],  
\end{align}
}we have for the solution $(U,V)$ of the corresponding initial value problem for \eqref{eq:PDE_model} that
\begin{equation}
    \left(U(\cdot,s),V(\cdot,s)\right) \in \Gamma^{\rightarrow}(\theta_1, \widehat{\theta}_2) \quad\text{for all}\quad s > s_0 \,.
\end{equation}
\end{res}

This result states that the front shape of initial conditions is quantitatively preserved by the PDE flow, but not necessarily qualitatively, since the $\Gamma^{\rightarrow}(\theta)$ is not forward invariant.

\begin{rmk}
 \label{rem:inTF_initialphase} In case the assumptions on the initial conditions are relaxed from having compact support to being sufficiently localised in space, a similar result can be formulated with the addition that the front will be formed after a short initial calibration phase that is of order $\mathcal{O}(\varepsilon^2)$ in $s$ (or equivalently order $\mathcal{O}(1)$ in $t$).
\end{rmk}

\paragraph{Proof strategy for mixed feedback $\alpha > 0$}\cite[Chapter 4]{LvV_thesis}
When $\alpha>0$, both boundedness between the background states (property \ref{prop:bounded_between_backgroundstates}) and uniform direction of propagation (property \ref{prop:propagation_direction}) are, in contrast to the positive feedback case, not forward invariant under the flow of \eqref{eq:PDE_model}. As a result, crucial building blocks that were used in the proof of Results \ref{res:TF_noSF} and \ref{res:TF_withSF} are missing.\\
However, the existence of an initialised travelling front -- i.e. a solution that maintains its front shape and has a uniform propagation direction for all $s >s_0$ -- can for $\alpha>0$ be shown using estimates on $\partial_s U$ based on the variation of constants formula applied to \eqref{eq:PDE_model_UtVt}. The assumption $\sup_{x \in \mathbb{R}} \frac{\mathrm{d}v_b^-}{\mathrm{d}x}(x;0) < \frac{\gamma}{\alpha}$ is crucial, implying that the integral $\int_{z(s) - \sqrt{\eps}}^{z(s)+\sqrt{\eps}} \partial_s U(x,s)\,\text{d}x$ is negative and uniformly bounded away from zero. Next, $\partial_s U$ is again written as the sum of three parts: a `fade out' part, which decreases exponentially in time; an `interface' part, whose properties are determined by the shape of the front interface, and a `tails' part, whose properties are determined by the asymptotic convergence to the background states. The assumption of compact support (Result \ref{res:inTF_compactsupport}) or the assumption of initial calibration (Result \ref{rem:inTF_initialphase}) can now be used to estimate the integral of $\partial_s U$ in an $\eps$-neighbourhood around the front position $z(s)$. Together with estimates on the `fade out' part and the `tails' part of $\partial_s U$, this provides a bound of the form $\partial_s U < - \frac{\lambda}{\eps}$ for some $\lambda>0$ for finite $s$. Next, the uniformity of the obtained estimates can be used to extend the $s$-interval of existence to all $s \geq s_0$, thereby proving the existence of initialised travelling front solutions to 
\eqref{eq:PDE_model_fast_reaction}.

\subsubsection{Entire travelling front solutions}\label{ss:TF_entire}
To make statements about the properties of travelling fronts with nonzero speed, it is necessary to define what we mean by "position":

\begin{deftn}\label{def:position}
    The position of an entire travelling front $(U_\text{TF},V_\text{TF})$ is defined as the unique function $z(s)$ for which $U_\text{TF}(z(s),s) = 0$ for all $s \in \mathbb{R}$.
\end{deftn}
The existence and uniqueness of $z(s)$ follow from the front properties stated in Definition \ref{def:front}. For stationary fronts, we have $z(s) \equiv x_0$, cf. Result \ref{res:stationaryfronts}. If no stationary fronts exists, the existence of uniformly travelling fronts is summarised in the following Result:

\begin{res}[\bf Entire travelling fronts]\label{res:TF_noSF}
 Let $\alpha \neq 0$, $f_1 = 0$ and \textcolor{black}{let} $f_2$ satisfy Assumption~\ref{ass:f}.
    \begin{enumerate}
        \item Let $\sup_{x \in \mathbb{R}} \alpha\frac{\mathrm{d}v_b^-}{\mathrm{d}x}(x;0) < \gamma$. Then there exists an $\eps_0>0$ such that for all $0<\eps\leq\eps_0$, there exists an entire travelling front. Its position satisfies $\frac{\mathrm{d}z}{\mathrm{d}s}>0$; moreover, $z(s) \to \pm \infty$ as $s \to \pm \infty$, and $\left\|\frac{\mathrm{d}z}{\mathrm{d}s}\right\|_\infty = \mathcal{O}(1)$ in $\eps$.
        \item Let $\inf_{x \in \mathbb{R}}\alpha \frac{\mathrm{d}v_b^-}{\mathrm{d}x}(x;0) > \gamma$. Then there exists an $\eps_0>0$ such that for all $0<\eps\leq\eps_0$, there exists an entire travelling front. Its position satisfies $\frac{\mathrm{d}z}{\mathrm{d}s}<0$; moreover, $z(s) \to \mp \infty$ as $s \to \pm \infty$, and $\left\|\frac{\mathrm{d}z}{\mathrm{d}s}\right\|_\infty = \mathcal{O}(1)$ in $\eps$.
    \end{enumerate}
\end{res}

We refrain from giving a full proof, which can be found in \cite{LvV_thesis}.

\paragraph{Sketch of the proof for Result~\ref{res:TF_noSF}}

\noindent We fix a non-negative compactly supported smooth even function $\Lambda$. The support (more precisely the complement of the zero level set) of $\Lambda$ is an interval centred around $0$, and its width is chosen such that it contains $u_h^{-1}([-\frac{1}{\sqrt{3}},\frac{1}{\sqrt{3}}])$, where $u_h$ is the Allen-Cahn front $u_h(\xi)= \tanh(\frac{1}{2}\sqrt{2} \xi)$. Given $z_0\in \mathbb{R}$, we construct an initial condition $(U_0,V_0)$ at position $z_0$ by enforcing that for some scalar $c\in \mathbb{R}$ the following relation holds
\begin{eqnarray}
(U_0^{(1)},V_0^{(1)})(x)= -c \Lambda\left(\frac{x-z_0}{\epsilon}\right) (U_0^\prime(x),V_0^\prime(x)) \, .\label{U01U0p}  
\end{eqnarray}	

\noindent  One can think of \eqref{U01U0p} as corresponding to a mixture of using a co-moving frame $\xi=\frac{x-(z_0+c(s-s_0) )}{\epsilon}$ near the interface and resting frame $\xi=\frac{x-z_0}{\epsilon}$ in the tails. The interface is cut off using the function $\Lambda$. The construction is instantaneous, in the sense that the position $z_0$ and the speed $c$ are fixed. Hence, we are interested in the 'germ' of stationary solutions of the PDE \eqref{eq:PDE_model} in terms of the (mixed co-moving/resting) coordinates $(\xi,s)$, i.e.  when $s$ is infinitesimally close to $s_0$, which gives (\ref{U01U0p}). More details can be found in  \cite[Chapter 4]{LvV_thesis}.\\ 
   
\noindent  Substituting \eqref{U01U0p} into the PDE gives rise to a system of nonautonomous ODEs.   In the system of ODEs $z_0\in \mathbb{R}$ appears as a parameter, and the aim is to find  $\varepsilon_0>0$ and a function $c=c(z_0,\varepsilon)$ (extending smoothly to $\varepsilon=0$), such that for $c=c(z_0,\varepsilon)$ and $\varepsilon \in (0,\varepsilon_0)$ the system of ODEs has an (appropriately defined) heteroclinic orbit. In this sense the construction parallels the construction of stationary fronts, where $c= c(z_0,\varepsilon)$ takes the role of $\gamma=\gamma_{SF}(z_0,\epsilon)$ (see Result~\ref{res:stationaryfronts}). However, a new technical challenge in the present case (which did not appear in the construction of stationary fronts) is that the system of ODEs is not a slow/fast system, due to the term $\Lambda(\frac{x-z}{\varepsilon})$, which either makes the vector field discontinuous at $\varepsilon=0$ or forces to introduce a variable $\zeta$ which in the fast scale satisfies $\zeta_\xi=1$ (in particular $\zeta$ is neither a fast nor a slow variable). Nevertheless, by making suitable adjustments to the (slow/fast) approach for stationary fronts, we can prove the existence of $\varepsilon_0>0$ and the function $c(z_0,\varepsilon)$ defined for $(z_0,\varepsilon)\in \mathbb{R}\times [0,\varepsilon_0]$.\\
    
\noindent For purposes of illustration, consider mixed feedback $\alpha>0$ and $\gamma$ as in part 1 of Result~\ref{res:TF_noSF} (i.e. the case of fronts moving with positive non-zero speed). The assumptions on $\Lambda$ and boundedness of the function $c(z_0,\varepsilon)$ implies that the resulting initial conditions $x\to (U_0,V_0)(x;z_0,\varepsilon)$ (extracted from the ODE heteroclinics) are elements of $\Gamma^\rightarrow(\theta)$. It can be shown that in fact the $m^+$ and $n^+$ components of $\theta$ are zero and that the   assumptions for Result~\ref{res:inTF_compactsupport} hold (in particular $U^{(1)}_0$ has compact support by construction).  Hence by Result~\ref{res:inTF_compactsupport} (for mixed feedback) there exist $\widehat{\theta}$ such that the solution $(U,V)(x,s;s_0,z_0,\varepsilon)$ of the IVP with initial condition $(U_0,V_0)(x;z_0,\varepsilon)$ at time $s_0$ satisfies  $(U,V)(x,s;s_0,z_0,\varepsilon)\in\Gamma^\rightarrow(\widehat{\theta}) $ for all $s\ge s_0$, which includes a uniform bound on the gradient of $(U,V)$. In addition, it can also be seen that the decay in the tails is uniform for all $s\ge s_0$, $s_0\in \mathbb{R}$ and $z_0\in \mathbb{R}$. The entire solution $(U,V)_{TF}$, is extracted from the family $\{(U,V)(x,s;s_0,z_0,\varepsilon): s_0\in \mathbb{R}, z_0\in \mathbb{R} \}$ using the Arzel\`{a}-Ascoli theorem along with the uniform estimates for $(U,V)$, where the uniformity for the length of the interface and the decay in the tails helps to deal with the lack of compactness of the domain $(x,s)\in \mathbb{R}^2$.


\paragraph{Stationary fronts as attractor of travelling fronts for $\alpha < 0$.}
For $\alpha < 0$ (positive feedback) one can establish that, if stationary fronts do exist, they act as attractors for travelling front in forward and backward time, as summarised in the following Result. \textcolor{black}{Its proof can be found in \cite[Chapter 4]{LvV_thesis}, where it is seen to follow naturally from the Theorem corresponding to Result \ref{res:TF_noSF}.}

\begin{res}\label{res:TF_withSF}
 Let $\alpha < 0$, $f_1 = 0$, \textcolor{black}{and let} $f_2$ satisfy Assumptions \ref{ass:f}. Let $(U_\mathrm{SF},V_\mathrm{SF})$ be a stationary front with location $x_0$, such that $\frac{\mathrm{d}^2 v_b^-}{\mathrm{d}x^2}(x_0;0) \neq 0$.
 \begin{enumerate}
     \item Suppose that $\frac{\mathrm{d}^2 v_b^-}{\mathrm{d}x^2}(x_0;0)>0$. Then there exists an $\eps_0>0$ such that for all $0<\eps\leq\eps_0$ there exists a pair of entire travelling fronts $(U_\mathrm{TF}^{\leftarrow/\rightarrow},V_\mathrm{TF}^{\leftarrow/\rightarrow})$ for which $\lim_{s \to -\infty} (U_\mathrm{TF}^{\leftarrow/\rightarrow},V_\mathrm{TF}^{\leftarrow/\rightarrow}) = (U_\mathrm{SF},V_\mathrm{SF})$. The front velocities satisfy $\frac{\mathrm{d}z^{\leftarrow}}{\mathrm{d}s} < 0$ and $\frac{\mathrm{d}z^{\rightarrow}}{\mathrm{d}s} > 0$; moreover, we have $u_b^-(x) \leq U^\rightarrow_\mathrm{TF}(x,s) < U_\mathrm{SF}(x) < U^\leftarrow_\mathrm{TF}(x,s) \leq u_b^+(x)$ and $v_b^-(x) \leq V^\rightarrow_\mathrm{TF}(x,s) < V_\mathrm{SF}(x) < V^\leftarrow_\mathrm{TF}(x,s) \leq v_b^+(x)$.
     \item Suppose that $\frac{\mathrm{d}^2 v_b^-}{\mathrm{d}x^2}(x_0;0)<0$. Then there exists an $\eps_0>0$ such that for all $0<\eps\leq\eps_0$ there exists a pair of entire travelling fronts $(U_\mathrm{TF}^{\leftarrow/\rightarrow},V_\mathrm{TF}^{\leftarrow/\rightarrow})$ for which $\lim_{s \to \infty} (U_\mathrm{TF}^{\leftarrow/\rightarrow},V_\mathrm{TF}^{\leftarrow/\rightarrow}) = (U_\mathrm{SF},V_\mathrm{SF})$. The front velocities satisfy $\frac{\mathrm{d}z^{\leftarrow}}{\mathrm{d}s} < 0$ and $\frac{\mathrm{d}z^{\rightarrow}}{\mathrm{d}s} > 0$; moreover, we have $u_b^-(x) \leq U^\leftarrow_\mathrm{TF}(x,s) < U_\mathrm{SF}(x) < U^\rightarrow_\mathrm{TF}(x,s) \leq u_b^+(x)$ and $v_b^-(x) \leq V^\leftarrow_\mathrm{TF}(x,s) < V_\mathrm{SF}(x) < V^\rightarrow_\mathrm{TF}(x,s) \leq v_b^+(x)$.
 \end{enumerate}
\end{res}

\newpage


\section{Delay-differential equation for front position dynamics}\label{s:delay_equation}

The results of the previous section provide general properties of the front position $z = z(s)$ from Definition~\ref{def:position}. In order to study the dynamics of $z(s)$ in more detail, one can derive an evolution equation for the front position. The outcome of this analysis, as carried out in \cite[Chapter 5]{LvV_thesis}, is formulated in the following Result:


\begin{res}[\bf Delay-differential equation]\label{res:delay_eq}
    Let $f_1 = 0, \eps_0>0$ and let $\left(U_\mathrm{TF}(x,s;\eps),V_\mathrm{TF}(x,s;\eps)\right)$ be a family of entire travelling front solutions (as given in Results \ref{res:TF_noSF}) parametrised by $\eps \in (0,\eps_0)$, with front position $z(s;\eps)$. The singular limit $z_R(s) := \lim_{\eps \downarrow 0} z(s;\eps)$ obeys the equation
    \begin{equation}\label{eq:delay_eq}
        \frac{\sqrt{2}}{3} \frac{\mathrm{d}z_R}{\mathrm{d}s} = \gamma + \alpha \left(\frac{\mathrm{d}v_b^+}{\mathrm{d}x}(z_R(s);0) + \hat{\tau} W[z_R]\textcolor{black}{(s)}\right),
    \end{equation}
    where $W[z]$ is the functional
    \begin{align}\label{eq:delayfunctional}
    W[z](s) := \frac{1}{2\sqrt{\pi\hat{\tau}}}
    \int_{-\infty}^s \int_{-\infty}^0 
        \frac{\mathrm{d}z}{\mathrm{d}s}(r)\left(1+f_2(z(r))\right)
        \,\textcolor{black}{\frac{1}{\sqrt{s-r}}}\,
        e^{x+\frac{r-s}{\hat{\tau}}}
        \left[
            e^{-\frac{\hat{\tau}}{4(s-r)}(x-(z(s)-z(r)))^2}
            +
            e^{-\frac{\hat{\tau}}{4(s-r)}(x+(z(s)-z(r)))^2}
        \right]
    \,\mathrm{d}x\, \mathrm{d}r.
\end{align}
\end{res}

Equation \eqref{eq:delay_eq} is derived rigorously in \cite[Chapter 5]{LvV_thesis} in the case $\alpha<0$, for which the existence of the $\eps$-family of entire front solutions is provided by Results \ref{res:TF_noSF} and \ref{res:TF_withSF}. For general $\alpha$, assuming the existence of a $\eps$-family of front solutions, equation \eqref{eq:delay_eq} can be formally derived using the approach sketched below. Numerical simulations suggest that \eqref{eq:delay_eq} indeed also holds for $\alpha > 0$.

\begin{rmk}\label{rmk:implicitdelay}
An implicit version of the above delay-differential equation is given by
\begin{align}
    & \alpha V_R(z_R(s),s) + \gamma = \frac{\sqrt{2}}{3} z_R'(s) \, ,\\
    & \hat{\tau} \partial_s V_R = \partial_x^2 V_R - V_R + (1+f_2(x))\mathrm{sign}(x-z_R(s)) \, .
\end{align}
In \cite{BC-BD.2020} a coupled algebraic-ODE system \textcolor{black}{was} formally derived to describe the dynamic position of travelling pulse solutions that bifurcated through a drift instability from a stationary one. Here we get for $\hat{\tau} > 0$ a PDE  which boils down to an ODE for $\hat{\tau} = 0$.    
\end{rmk}

\subsection{Formal derivation}
Before showing the formal derivation, we provide a short overview. To derive \eqref{eq:delay_eq}, we mimic the approach for stationary fronts from Section~\ref{ss:stationaryfronts}, writing \eqref{eq:PDE_model} \textcolor{black}{as} a "time-dependent spatial dynamics system" (with $U_{TF},V_{TF}$ as in Result~\ref{res:delay_eq})
 \begin{subequations}\label{eq:delayeq_ODEsystem}
     \begin{align}
     \eps u_x &= p,\\
     \eps p_x &= -u + u^3 + \eps(\alpha v+\gamma) \textcolor{black}{+ \eps^2 \partial_s U_{TF}(x,s;\varepsilon)},\\
     v_x &= q,\\
     q_x &= v - (1+f_2(x))u \textcolor{black}{+ \hat{\tau} \partial_s V_{TF}(x,s;\varepsilon)}.
     \end{align}
 \end{subequations}

\noindent We are interested in the singular limit $\varepsilon=0$, which gives, with $z_R$ as in Result~\ref{res:delay_eq}

\begin{subequations}\label{eq:delayeq_ODEsystemeps0}
     \begin{align}
   v_x &= q,\\
     q_x &= v - (1+f_2(x))\textnormal{sign}(x-z_R(s)) \textcolor{black}{+ \hat{\tau} \partial_s V_R(x,s)}.
     \end{align}
 \end{subequations}

\noindent Treating $\hat{\tau} \partial_s V_R$ as an inhomogeneity of the $(v,q)$-system, we find (cf. \eqref{eq:solutionoperator_Veq}) to leading order
\begin{equation}
    V_R(x,s) = \left\{\begin{array}{rcl} v_b^-(x;0) \textcolor{black}{-} G\left(\hat{\tau} \partial_s V_R(\cdot,s)\right) -(v_b^-(z_R(s);0) \textcolor{black}{+} \frac{\text{d} v_b^-}{\text{d}x}(z_R(s);0)) e^{x-z_R(s)}& \text{for} & x < z_R(s) \\[.2cm]  v_b^+(x;0) \textcolor{black}{-} G\left(\hat{\tau} \partial_s V_R(\cdot,s)\right) -(v_b^+(z_R(s);0) \textcolor{black}{+} \frac{\text{d} v_b^-}{\text{d}x}(z_R(s);0)) e^{-x+z_R(s)} & \text{for} & x > z_R(s) \end{array}\right.
\end{equation}
Taking the partial derivative w.r.t. $s$ and equating $\partial_s v_R(x,s) = \partial_s V_R(x,s)$ yields the following PDE for $w(x,s) := \textcolor{black}{-}G\left(\partial_s V_R(\cdot,s)\right)$:
\begin{equation}\label{eq:delayeq_PDEw}
    \hat{\tau}\partial_s w = \partial^2_x w - w + z_R'(s)(1+f_2(z_R(s))) e^{-|x-z_R(s)|}
\end{equation}

 \noindent By analogy with the constant coefficient case 
 $f_2=0$ it seems plausible that to leading order in $\eps$ the following Melnikov condition should hold

\begin{equation}
 \frac{\sqrt{2}}{3} \frac{\text{d} z_R}{\text{d}s} = \alpha V_R(z_R(s),s)+ \gamma.   
\end{equation}
Using the explicit form of the heat semigroup on \eqref{eq:delayeq_PDEw} yields Result \ref{res:delay_eq}.\\

\noindent We will discuss the formal derivation in more detail. For convenience, we drop the index $R$ in the notation $V_R$ and $z_R$, which will not result in confusion as this section involves only the singular limit $\varepsilon=0$. We will also use the notation $q_b^\pm:=\frac{\mathrm{d}v_b^\pm}{\mathrm{d}x}$

\paragraph{Preliminaries}\label{ss:DDE_prelim}
The starting point is to borrow intuition from the constant coefficient case when $f_2=0$. In that case, there exist for $0<\varepsilon \ll 1$ travelling fronts which travel with constant speed $c(\varepsilon)$, which in the singular limit $\varepsilon \downarrow 0$ satisfies \eqref{eq:existence_const_coeff}, which we recall here for clarity
\begin{eqnarray}
\alpha v^*(c) + \gamma = \frac{\sqrt{2}}{3} c,\quad v^*(c)=\frac{\hat{\tau}c }{\sqrt{4+c^2\hat{\tau}^2}} \label{numerics_Melnikov}
\end{eqnarray}

\noindent where $v^*(c)$ is the value of $V$ at the interface of the step function $U= \textnormal{sign}(x-z(s))$, and $z(s)=z_0+cs$ is the position of the interface at time $s$. By analogy, it is plausible that when the heterogeneity $f_2$ is introduced, the following existence condition should hold
 (justifying this rigorously requires much more work, in particular proving existence of travelling fronts $(U,V)$ for $0<\varepsilon \ll 1$ and studying the  limit behaviour $(U,V)$ as well as the gradient $(U_x,U_s,V_x,V_s)$ as $\varepsilon \downarrow 0$,  and is one of the main topics of \cite{LvV_thesis}.) 
\begin{eqnarray}
	\alpha V(z(s),s) + \gamma = \frac{\sqrt{2}}{3}z^\prime(s)  
	.\label{numerics_av2}
\end{eqnarray}
To compute $V(z(s),s)$ we note that $V(x,s)$ should be a front solution of the inhomogeneous scalar PDE 
\begin{align}
	\hat{\tau}V_s =V_{xx}-V+ (1+ f_2(x)) \textnormal{sign}(x-z(s)) \label{PDEVseps0}
\end{align}
 where the inhomogeneity has a discontinuity at $x=z(s)$. The first intuition would be to write down the  variation of constants formula
\begin{align}
&	V(x,s):= \int_{-\infty}^s \int_{-\infty}^{\infty} \frac{1}{\hat{\tau}}(1+ f_2(x-y)) \text{sign} (x-y-z(u))  \frac{\sqrt{\hat{\tau}}}{\sqrt{4\pi}\sqrt{s-u}}\exp(-\frac{y^2\hat{\tau}}{4(s-u)} )\exp(-\frac{s-u}{\hat{\tau}})  dy du,
	\end{align}
however, it turns out to  be difficult  to distil a meaningful expression for $V(z(s),s)$ in a direct way (i.e. by manipulating the integrals using change of variables, partial integration etc).
	
\paragraph{GSPT approach to the time-dependent spatial dynamics system \eqref{eq:delayeq_ODEsystem}}
We will look from a geometric point at the situation and derive a relation (see \eqref{numerics_Gv2}) between $V(z(s),s)$ and $G(V_s(\cdot,s))(z(s))$, where  (see Result~\ref{res:backgroundstates})
		\begin{eqnarray}
		G(f)(x): = -\frac{1}{2}\left( \int_x^\infty e^{-y}f(y)dy   \right) e^x - \frac{1}{2}\left( \int_{-\infty}^x e^{y} f(y)dy     \right) e^{-x}. \nonumber 
		\end{eqnarray}
is the solution operator corresponding to the ODE $v^{\prime\prime}-v=f$. To this end we consider for fixed $s\in \mathbb{R}$ the reduced slow systems given by
\begin{align}
	\begin{split} 
	&v^{\prime}   =  q \\
	&q^{\prime} = v - (1+f_2(x))u  + \hat{\tau} V_s(x,s) \\
	& u \in \{-1,1\}
	\end{split}  \label{numerics_redslow},
	\end{align}
where we view $V_s(\cdot,s)$ \textcolor{black}{as} an inhomogeneity (and for the moment ignoring the relation between $V_s$ and $V$). If the term $V_s(x,s)$ would not be present, then (under assumptions A1 and A2 on $f_2$) the system  \eqref{numerics_redslow}  with $u=-1$ has a unique bounded negative solution (i.e. singular limit background state) $(v_b^-,q_b^-)(x,0)$, $x\in \mathbb{R}$  and a unique bounded positive solution $(v_b^+,q_b^+)(x,0)$, $x\in \mathbb{R}$. The solution $V$ of \eqref{PDEVseps0} which we would like to construct  satisfy the following properties. 
	\begin{enumerate} 
	\item (connect negative background state in backward time) $\{(V(x,s),V_x(x,s)), x\le z(s)\}$ which is   a (backward) orbit of \eqref{numerics_redslow} for $u=-1$, should connect to $(v_b^-,q_b^-)(x,0)$ as $x\to -\infty$.
	\item (connect positive background state in forward time) Likewise
   $\{(V(x,s),V_x(x,s))$, $x\ge z(s)\}$ is an orbit of \eqref{numerics_redslow} for $u=1$ which should connect to $v_b^+(x,0)$ as $x\to \infty$. 
   \item (heteroclinic)  These backward and forward orbits should match at $x=z(s)$ 
   \item (Melnikov condition) $V(z(s),s)$ should satisfy the Melnikov condition \eqref{numerics_av2}
   \end{enumerate} 
We will interpret the effect of the inhomogeneity $x\to V_s(x,s)$ in \eqref{numerics_redslow} as a distortion of the background states $(v_b^\pm (x;0))$ near the interface $x=z(s)$ such that the construction described can be carried out for all $s\in \mathbb{R}$. Indeed, if the inhomogeneity $V_s(\cdot,s)$ in \eqref{numerics_redslow} would not be present, then geometrically, the conditions imply that  the three lines $(v_b^-,q_b^-)(z(s);0)+ \mathbb{R}(1,1)$ and $-(v_b^+,q_b^+)(z(s),0)+ \mathbb{R}(1,-1)$, and $\alpha v + \gamma = \frac{\sqrt{2}}{3}z^\prime(s)$ in the $(v,q)$ plane have to intersect for all $s\in \mathbb{R}$. This situation is very rigid, e.g. even if we consider just a single value of $s$, we still need three lines to intersect, which we can only expect to happen at a small number of positions $z(s)=s$. In particular setting $z^\prime(s)=0$ and $z(s)$ is equal to (a small number of) specific positions, leads to the construction of stationary fronts. \\
    
   \noindent In view of the preceding discussion, define
	\begin{align}
	\begin{split} 
	&v_{b,\textnormal{dist.}}^{\pm }(x;s): =   v_b^{\pm}(x,0) + \widetilde{v}(x;s) \\
	&\widetilde{v}(x;s):=\hat{\tau} G(V_s(.,s))(x)
	\end{split}
	\label{bsoldist}
	\end{align}
	here we interpret $v_{b,\textnormal{dist}}^{\pm }$ as distorted background states. Note that $v_b^{\pm}(x,0)  = - (\pm) G ( 1+f_2(x) )$. Define $\tilde{q}(x;s)= \frac{\mathrm{d}}{\mathrm{d}x}\tilde{v}(x;s)$ and $q_{b,\textnormal{dist}}^{\pm }(x;s)= \frac{\mathrm{d}}{\mathrm{d}x}\tilde{v}(x;s)$. The effect of the distortion is to make the following lines in the $(v,q)$-plane intersect
	\begin{itemize} 
	 \item The vertical line $\alpha v + \gamma = \frac{\sqrt{2}}{3}z^\prime(s)$. Let
	 \begin{align}
	 v^*(s) :=   \frac{1}{\alpha}\left( \frac{\sqrt{2}}{3}z^\prime(s) - \gamma\right).\nonumber 
	 \end{align}
	\item The line
    \begin{align}\label{l1} 
	   \left\{ (v,q) \in \mathbb{R}^2 ~|~ (v,q)= -(v_b^-,q_b^-)(z(s);0) + (\tilde{v}(z(s);s), \tilde{q}(z(s);s)) + \mathbb{R}(1,1) \right\},  
	\end{align}
	which if we make \eqref{numerics_redslow}, $u=-1$ autonomous by setting $\chi_x=1$, we may interpret as the unstable manifold of the distorted background state $(v_{b,\textnormal{dist}}^{-}, q_{b,\textnormal{dist}}^{- })(x;s)$ viewed in the slice $\{\chi=z(s)\}$.
	\item  The line 
    \begin{align}\label{l2} 
	   \left\{ (v,q) \in \mathbb{R}^2 ~|~ (v,q)= (v_b^-,q_b^-)(z(s),0) + (\tilde{v}(z(s);s), \tilde{q}(z(s);s)) + \mathbb{R}(1,-1)  \right\},  
	\end{align}	
	which if we make \eqref{numerics_redslow}, $u=+1$ autonomous by setting $\chi_x=1$, we may interpret as the stable manifold of the distorted background state $(v_{b,\textnormal{dist}}^{+}, q_{b,\textnormal{dist}}^{+ })(x;s)$ viewed in the slice $\{\chi=z(s)\}$.
	\end{itemize} 
	To see when the lines \eqref{l1} and \eqref{l2} intersect at $v=v^*$ we consider the following lemma. 
    
	\begin{lem}\label{abcd} 
		Let $a,b,c,d\in \mathbb{R}$ and let $\ell_i$, $i=1,2,3$ be $3$ lines of the form $\ell_1 : (v,q)= [(a,b)+(c,d)]  +\mathbb{R}(1,1)$, $\ell_2 :(v,q) =[-(a,b)+ (c,d)] + \mathbb{R}( 1,-1)$ and $\ell_3:v=v^*$. Then $\ell_i$ intersect in a single point if and only if $c=v^* + b$.
	\end{lem}
	\textbf{Proof:} We need to solve the (overdetermined) system of equations $a+ c+ t = v^*$, $-a+c+s= v^*$ and $b+d + t = -b + d -s$ for $s,t$. From the first two equations show $t= v^*-a-c$ and $s= v^*+a-c$. Substituting in the third equations gives the compatibility condition $c=  v^*+b$. $ \blacksquare$\\
	
	\noindent Applying Lemma~\ref{abcd} to (\ref{l1}) and (\ref{l2}) (and recalling $(v_b^+,q_b^+)=-(v_b^-,q_b^-)$ ) gives $\tilde{v}(z(s);s)   =  v^*(s) +q_b^-(z(s);0) $. Recalling $\tilde{v}(z(s);s)= \textcolor{black}{-}\hat{\tau}G(V_s(\cdot,s))(z(s))$ (see \eqref{bsoldist}), and using that $V(z(s),s)$ should be equal to $v^*(s)$ we obtain the relation
	\begin{eqnarray}
	V(z(s),s)=  \tilde{v}(z(s);s)   - q_b^-(z(s),0)= \textcolor{black}{-}\hat{\tau}G(V_s(\cdot,s))(z(s))   - q_b^-(z(s),0). \label{numerics_Gv2} 
	\end{eqnarray}
    \noindent In addition we have
	\begin{eqnarray}
	V(x,s) = v_{b,\textnormal{dist.}}^{-}(x;s) +  A \exp(x)  \text{ if $x\le z(s)$} \nonumber \\
	V(x,s) = v_{b,\textnormal{dist.}}^{+}(x;s) + B \exp(-x) \text{ if $x\ge z(s)$},  \nonumber 
	\end{eqnarray}	
	\noindent for constants $A$, $B$.   
	Using that $v_{b,\textnormal{dist}}^{\pm }(z(s);s)=  v_b^{\pm}(z(s);0) + V(z(s);s) + q_b^{-}(z(s);0)$  (in view of \ref{numerics_Gv2}), we obtain 
\begin{eqnarray} 
	A=-(v_b^-(z(s);0)+ q_b^-(z(s);0) \exp(-z(s))\nonumber \\
	B= -(v_b^+(z(s);0)+ q_b^-(z(s);0) \exp(z(s))\nonumber 
	\end{eqnarray} 
	\noindent In particular note that $q_b^-(z(s);0)$ appears in both equations (so it is not $q_b^-$ in one equation and $q_b^+$ in the other). In summary we have found
	\begin{eqnarray}
	V(x,s)   = \left\{
	\begin{array}{c}
	v_b^-(x;0)  \textcolor{black}{-}\hat{\tau} G(V_s(\cdot,s))(x) 	-(v_b^-(z(s);0)+ q_b^-(z(s);0))\exp(x-z(s))\text{ if $x\le z(s)$},\\
	v_b^+(x;0)  \textcolor{black}{-} \hat{\tau} G(V_s(\cdot,s))(x)  -(v_b^+(z(s);0)+ q_b^-(z(s);0))\exp(-(x-z(s))) \text{ if $x\ge  z(s)$}
	\end{array} \right. \label{vxs} 
	\end{eqnarray}
    
\paragraph{Deriving an expression for $W[z_0]$}
So far we viewed $s$ as fixed and looked at the profile $x\to V(x,s)$. We will now look at $V$ as a function of two variables $(x,s)$.
	 Set $\omega(x,s) := \textcolor{black}{-}G(V_s(\cdot,s))(x)$. We differentiate \eqref{vxs} w.r.t $s$. On the left hand side of \eqref{vxs} we get $V_s(x,s)= \omega_{xx}-\omega$.\\
	
\noindent  On the right hand side we compute 
	\begin{align}
	\begin{split} 
	& 
	\partial_s \left[  (v_b^-(z(s);0)+ q_b^-(z(s);0))\exp(x-z(s))\right]   =\\
	& z^\prime(s) (q_b^-(z(s);0) + v_b^-(z(s);0) + (1+ f_2(z(s))) )\exp(x-z(s)) \\
	&- z^\prime(s)  (v_b^-(z(s);0)+ q_b^-(z(s);0))\exp(x-z(s)) \\
	&= z^\prime(s) ( 1+f_2(z(s)))\exp(x-z(s) )
	\end{split}    \nonumber
	\end{align}
	\noindent Likewise we compute using $(v_b^+,q_b^+)=-(v_b^-,q_b^-)$
	\begin{align}
	\begin{split} 
	& 
	\partial_s\left[  (v_b^+(z(s);0)+ q_b^-(z(s);0))\exp(-(x-z(s)))\right]   =\\
	& z^\prime(s) \bigl(q_b^+(z(s);0) + v_b^+(z(s);0) - (1+ f_2(z(s))) \bigr)\exp(-(x-z(s))) \\
	&+ z^\prime(s)  (v_b^+(z(s);0)+ q_b^-(z(s);0))\exp(-(x-z(s))) \\
	&= z^\prime(s) ( 1+ f_2(z(s)))\exp(-(x-z(s)) )
	\end{split}    \nonumber
	\end{align} 
 We have
 \begin{eqnarray}
	\partial^2_x \omega - \omega = 
	 \hat{\tau} \partial_s \omega +  h(x,s), \nonumber
	\end{eqnarray} 	
where 		
	\begin{eqnarray}
		h(x,s)=  - z^\prime(s)(1 + f_2(z(s))) \exp(-|x-z(s)|) \nonumber 
	\end{eqnarray}	
Applying the Fourier transform (with respect to the variable $x$) and rearranging
\begin{eqnarray}
	\hat{\tau}  \partial_s \textcolor{black}{\hat{\omega}} +(k^2+1)\textcolor{black}{\hat{\omega}} = -\hat{h}(k,s) \nonumber
	\end{eqnarray}	
For each fixed $k$ this is a linear first order ODE. Denoting $\textcolor{black}{\hat{\omega}_0(k):=\hat{\omega}(k,0)}$ we have 
\begin{eqnarray}
	\textcolor{black}{\hat{\omega}}(k,s) = \exp\left(-  \frac{k^2+1}{\hat{\tau}}s   \right) 
	\left[\textcolor{black}{\hat{\omega}}_0(k) - \int_0^s \frac{1}{\hat{\tau}}\hat{h}(k,u)\exp \left(  \frac{k^2+1}{\hat{\tau}}u  \right)du \right] \nonumber
\end{eqnarray}	
We determine $\textcolor{black}{\hat{\omega}}_0(k)$ by imposing that the solution $s\to \textcolor{black}{\hat{\omega}}(k,s)$ should be bounded. Assuming $\hat{\tau}>0$ we can write
\begin{eqnarray}
	\hat{\omega}(k,s) = \exp\left(-  \frac{k^2+1}{\hat{\tau}}s   \right)
	\left[ - \int_{-\infty}^s \frac{1}{\hat{\tau}}\hat{h}(k,u)\exp \left(  \frac{k^2+1}{\hat{\tau}}u\right)du  \right] \nonumber
	\end{eqnarray}	
Using  $\mathcal{F}(x\to \exp(-|x-z(u)|)) =\textcolor{black}{e^{-ik z(u)}}\mathcal{F}(x\to \exp(-|x|))$ we obtain
\begin{align}
	&\hat{\omega}(k,s) = \mathcal{F}( \exp(-|.|) )\int_{-\infty}^{s}\frac{1}{\hat{\tau}} z^\prime(u)(1+ f_2(z(u))) \exp \left( \frac{k^2+1}{\hat{\tau}}(u-s ) - ik (z(u))  \right) du \nonumber
\end{align} 	
We now examine the inverse Fourier transform of the integral in the second line. For  $u<s$ we have 
\begin{eqnarray}
	\frac{1}{\sqrt{2\pi}}\int_{-\infty}^{\infty} e^{ikx} \exp \left( \frac{k^2+1}{\hat{\tau}}(u-s ) - ik (z(u))  \right) dk 
	= \frac{ \sqrt{\hat{\tau}}  }{ \sqrt{2}\sqrt{ s-u } } \exp\left(  \frac{ (z(u) -x)^2\hat{\tau}  }{4  (u-s) }  \right) \exp( \frac{u-s}{\hat{\tau}} )    \nonumber
\end{eqnarray}		
We obtain 
\begin{eqnarray}
			\omega(y,s) = \frac{1}{\sqrt{2\pi}}\int_{-\infty}^{s}\int_{-\infty}^{0}\frac{1}{\hat{\tau}} z^\prime(u)(1+ f_2(z(u)))\frac{ \sqrt{\hat{\tau}}  }{\sqrt{2} \sqrt{ s-u } }\exp( \frac{u-s}{\hat{\tau}} )    \nonumber\\
			e^x\left[\exp\left(  \frac{ -(-y+x+z(u) )^2\hat{\tau}  }{4  (s-u) }  \right) + \exp\left(  \frac{ -(-y-x+z(u))^2\hat{\tau}  }{4  (s-u) }   \right) \right]dx du    \nonumber
\end{eqnarray}
In particular,
		\begin{align}
		\begin{split} 
		\omega(z(s),s) =& \frac{1}{\sqrt{2\pi}}\int_{-\infty}^{s}\int_{-\infty}^{0}\frac{1}{\hat{\tau}} z^\prime(u)(1+ f_2(z(u)))\frac{ \sqrt{\hat{\tau}}  }{\sqrt{2} \sqrt{ s-u } }\exp( \frac{u-s}{\hat{\tau}} )\\
		&e^x\left[\exp\left(  \frac{ -(x-[z(s)-z(u)] )^2\hat{\tau}  }{4  (s-u) }  \right) + \exp\left(  \frac{ -(-x- [z(s)-z(u)] )^2\hat{\tau}  }{4  (s-u) }   \right) \right]dx du
		\end{split} \label{omegaz}
		\end{align}
Recalling $\omega: = \textcolor{black}{-}G(V_s) $, this completes the (formal) derivation of the delay differential equation, in view of \eqref{numerics_av2},  \eqref{numerics_Gv2}.


\subsection{Numerical algorithm for the delay-differential equation}
\noindent Let  $s_0$ be an initial time and let  $z_0\in \text{PC}^1((-\infty, s_0])$ be a function 
specifying a 'history'. {  We use the notation $z_0\in \mathrm{PC}^1((-\infty,s_0])$ for $z_0$ continuous and piecewise $C^1$ on $(-\infty,s_0]$
(with one-sided derivatives at the breakpoints).} Given $a\in \mathbb{R}$, we piecewise linearly extend $z_0(.)$ to a function $z[a](s)$, by setting  

	\begin{eqnarray}
	z[a](s)= z(s_0)+ (z^\prime(s_0)+a)(s-s_0) \text{   for $s\ge s_0$}.  \label{numerics_ca2} 
	\end{eqnarray}

\noindent Consider the delay functional $W[z]$ given by \eqref{eq:delayfunctional} and
the error term
	
	\begin{eqnarray}
	e(s,a):=\frac{\sqrt{2}}{3} z[a]^\prime(s)-\alpha( \hat{\tau}W[ z[a] ](s) - q_b^-(z[a](s));0) - \gamma, \nonumber 
	\end{eqnarray}

\noindent i.e. the extended function	$z[a](s)$ satisfies the delay differential equation at time $s$ if and only if $e(s,a)=0$.\\

\subsubsection{Algorithm 1}\label{s:delay_equation1}

\noindent \textbf{Input}: An initial time $s_0$ and a function  $z_0\in \text{PC}^1((-\infty, s_0])$ specifying a 'history',  such that $z_0(s)$ satisfies the delay differential equation at $s=s_0$. We do not require that $z_0(s)$ satisfies the delay differential equation at $s<s_0$ (from a practical point of view the algorithm even seems to work if $z_0(s_0)$ does not satisfy the delay differential equation,  and quickly converges in forward time $s$ to an actual solution) 
 Let $T>0$ and  consider the interval $[s_0,s_0+T]$. Let $N$ be an integer such that $h =\frac{T}{N}<h_0$  Denote $s_i:=s_0 + ih $, for $i=0,\ldots, N$.\\

 \noindent Set $\widehat{z}_0:=z_0(.)$.  To construct a numeric solution of the delay differential equation on $[s_0,s_0+T]$ we 
   proceed along the following step, for $i=1,\ldots N$

 \begin{enumerate}
 	\item Assuming $\widehat{z}_{i-1}\in \text{PC}^1((-\infty,s_{i-1}])$ has been defined on $(-\infty, s_{i-1}]$, consider the extended function $\widehat{z}_{i-1}[a]$. Consider the error term 
 	\begin{align}
 	e_i(s,a):= \begin{array}{l} \frac{\sqrt{2}}{3} \widehat{z}_{i-1}[a]^\prime(s)-\alpha( \hat{\tau}W[ \widehat{z}_{i-1}[a] ](s_i) - q_b^-(\widehat{z}_{i-1}[a](s_i));0) - \gamma.\end{array} \label{ei}  
 	\end{align}
 
 	\item We solve the equation $e_i(s_i,a)=0$, for $a$ as follows. Existence and uniqueness  of the solution $a=a_i^*$ is justified in \cite{LvV_thesis} for $0<h\ll 1$
 	(under suitable conditions on the history $\widehat{z}_0(.)$ involving Lipschitz continuity, see \cite{LvV_thesis}).

 	\item We  define $\widehat{z}_{i}\in \text{PC}^1((-\infty,s_{i}])$ by 
 	
 	\begin{eqnarray}
 	\widehat{z}_i := \widehat{z}_{i-1}[a_i^*]|_{(-\infty,s_i]} 
    \nonumber 
 	\end{eqnarray}

 	\noindent By construction we have $e_i(s_i,a_i^*)=0$, which means that $\widehat{z}_i(s)$ satisfies the delay differential equation at $s=s_i$. 	
 \end{enumerate}

 \noindent \textbf{Output}. Set $\widehat{z}:= \widehat{z}_N$. 
 Then (since $e_i(s_i,a_i^*)=0$ for all $i$ by construction),
 $\widehat{z}(s)$ solves the delay differential equation at $s=s_i$, for $i=0,\ldots ,N$.\\

\noindent In order to implement the algorithm one needs a way to efficiently compute  numerically the double integral involved in the delay differential equation. To this end we use the following probabilistic interpretation

\begin{eqnarray}
W[z](s)= 	 \mathbb{E}_{X\sim \text{Exp}(1) \atop R|X \sim \text{Levy}\left(0,  \frac{\hat{\tau}X^2}{2}\right) } \left[
\begin{array}{c} 
\frac{1}{\hat{\tau}X} z^\prime(s-R)(1+ f_2(z(s-R)))R\\ 	\exp\left(   -\left[\frac{1}{\hat{\tau}}  +\frac{\hat{\tau}}{4}  \left(\frac{z(s)-z(s-R)}{R}\right)^2\right]   R    \right)\\
2 \cosh\left(   \frac{1}{2} \hat{\tau}\left(\frac{z(s)-z(s-R)}{R}\right) X\right)
\end{array}
\right] .
   \nonumber
\end{eqnarray}
\noindent 
The probabilistic interpretation can be used to compute the double integral by means of using Monte Carlo integration. The main parameters determining the precision of the algorithm are the number of Monte Carlo samples $M$ per time steps, and the length of each time step $h$ (i.e. the temporal discretization) or equivalently, the temporal resolution $R=\lceil\frac{1}{h}\rceil$, i.e. the number of steps per unit time. Note that the computation time is $\mathcal{O}(MR)$.\\

\subsubsection{Algorithm 2}\label{s:delay_equation2}
A slight variation on Algorithm 1, which avoids Monte Carlo integration is based on the implicit form \eqref{rmk:implicitdelay} of delay-differential equation.
The idea is to numerically simulate $V_R$ alongside $z_R$ and involve the estimate of $V_R$ in the evaluation of the error term $e_i(s;a)$ (as in \eqref{ei}).\\

\noindent We start the simulation with an estimate of $V_R(.,s_0)$, 
which in principle can be obtained from the history $z(s)$, $s\le s_0$ via \eqref{vxs} and \eqref{omegaz} for $\omega= \mathcal{G}(\partial_sV_{R})$. An estimate for $V_R(.,s_0)$ can also be clear from the context (see e.g. example~0 in Section~\ref{ss:examples}). In addition to $\widehat{z}_0$ as described in the initialization step $i=0$ of Algorithm 1,  
set $\widehat{V}_{R,0}:= V_R(.,s_0)$. 
As in Algorithm 1, we proceed by iteration for $i=1,\ldots ,N$. In the $i$-th step of the simulation we consider the error
\begin{align}
e_i(s_i,a)= \frac{\sqrt{2}}{3} (\widehat{z}_{i-1}[a])^\prime (s) - \alpha   V (  \widehat{z}_{i-1}[a](s_i),s_i ;a) -\gamma\nonumber 
\end{align}
\noindent where  $V ( \cdot ;a)$ is obtained by numerically solving the initial value problem (with $\widehat{V}_{i-1}(\cdot;s_{i-1})$  defined in the previous iteration step $i-1$ )
\begin{align}
\hat{\tau} \partial_s V = \partial_x^2 V - V + (1+f_2(x))\mathrm{sign}(x-\widehat{z}_{i-1}[a] (s)), \quad V(\cdot,s_{i-1};a) = \widehat{V}_{i-1}(\cdot,s_{i-1})\nonumber 
\end{align}
\noindent for $s\in [s_{i-1},s_i]$. Solve $e(s_i,a)=0$ for $a$, giving $a=a_i^*$ and set 
\begin{align} 
& \left\{ \begin{array}{ll}
(\widehat{V}_i (\cdot,s), \widehat{z}_i(s)):= (\widehat{V}_{i-1}(\cdot,s), \widehat{z}_{i-1}(s)) & \text{ for $s\le s_{i-1}$}\\ (\widehat{V}_i(\cdot,s), \widehat{z}_i(s)) := (\widehat{V}(\cdot,s;a_i^*), \widehat{z}_{i-1}[a^*_i](s)) & \text{ for $s\in [s_{i-1},s_i]$} \nonumber 
\end{array}
\right.
\end{align}

{ 
\subsubsection{Algorithm 1 vs Algorithm 2} \noindent It is good to keep both Algorithm~1 and Algorithm~2 because they provide two genuinely different numerical routes to the same reduced front-position law, and that difference is useful rather than redundant. Algorithm~1 works directly with the explicit delay functional $W[z]$ and typically evaluates the defining double integral via its probabilistic (Monte--Carlo) interpretation, whereas Algorithm~2 uses the implicit formulation by evolving the auxiliary field $V_R$ through a PDE and coupling it back into the update of the position. Since these approaches have different sources of numerical error---Algorithm~1 is dominated by sampling variance (and time-stepping error), while Algorithm~2 is dominated by deterministic discretization and solver choices---keeping both provides an internal consistency check: when they agree (as in Example~0), this strongly supports both the reduction and the implementation. When they disagree (as in Example~3), the discrepancy is equally informative because it highlights parameter regimes in which one formulation or its discretization is more delicate, thereby aiding diagnosis of instability, stiffness, boundary effects, or insufficient resolution. Finally, the two methods offer complementary practical trade-offs: Algorithm~1 avoids building a full PDE solver and is convenient when one only needs $z(s)$, while Algorithm~2 can be efficient in some regimes, avoids Monte--Carlo sampling entirely, and simultaneously produces $V_R(\cdot,s)$, which is useful for interpreting and validating the front dynamics.

}

\subsection{Numerical exploration of front dynamics}\label{ss:examples}

{ 
We can now compare the dynamics predicted by the delay equation with direct simulations of the full PDE and explore several qualitatively different regimes. Example 0 considers a front in a heterogeneous medium without stationary fronts and shows how the heterogeneity modulates the front speed while preserving its direction of motion. Example 1 treats the case $\alpha<0$ with two stationary fronts and illustrates how travelling fronts connect and interact with these equilibria in agreement with Results~\ref{res:TF_noSF} and~\ref{res:TF_withSF}. Example 2 focuses on the homogeneous case $f_2\equiv 0$ with mixed feedback $\alpha>0$ and demonstrates direction-reversing fronts that connect different constant-speed travelling waves. Finally, Example 3 introduces a strong localized heterogeneity (outside the rigorous assumptions) that acts as a pair of “direction switches” and leads to numerically observed periodic front motion, providing an exploratory test of the delay equation beyond the theoretical setting. The corresponding code is avaiable at \cite{laraleiden2026code}.
}

\medskip

\paragraph{Example 0: Travelling front moving with a time-dependent speed $z^\prime(s)$ in the absence of stationary fronts.}
We consider a compactly supported $f_2$. The front enters and passes through the support of $f_2$. We set $\alpha=0.5$, $\gamma=0.2$, $\hat{\tau}=1$.
We define $f_2(x)$ as a piecewise polynomial { \sout{defined on $[-3,15]$}}. To this end, we sample the function
\begin{align}
(y-3)\sin\!\bigl(1+14(y-3)^3\bigr)\exp\!\bigl(-1.2|y-3|\bigr) \nonumber
\end{align}
\noindent for the $y$-values $y=3,4,5,6,7,8,9$ and then use MATLAB's \texttt{makima} interpolation to construct a piecewise polynomial defined on $[3,10]$ satisfying $f_2(3)=f_2(10)=0$. The sampled values $y=4,5,6,7,8,9$ (approximately) specify the local minima and maxima of the resulting $f_2$. For $x\le 3$ and $x\ge 10$ we set $f_2(x)=0$.

\begin{figure}[H]
    \centering
    	\scalebox{.6}{\includegraphics[trim = 2cm 9cm 4cm 9cm, clip,width=0.8\textwidth]{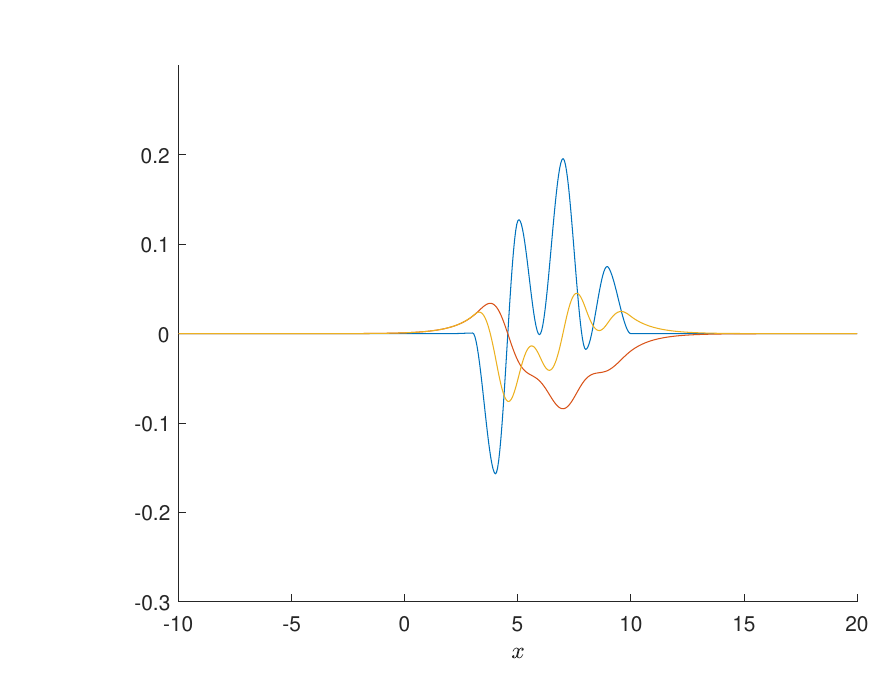}}
	\caption{ \small  $f_2(x)$ (in black), $v_b^-(x;0)+1$ (orange) and $q_b^-(x;0)$ (yellow). The shape of $f_2$ is chosen such that it involves multiple local minima and maxima. The corresponding $q_b^-(x;0)$ also fluctuates and we expect the front velocity to change over time as well, in an interesting way. The parameter $\gamma$ is chosen to avoid the presence of stationary fronts, i.e.\ intersections of the graph of $q_b^-(\cdot;0)$ with the horizontal line at height $\frac{\gamma}{\alpha}=0.4$.     
    }	
	\label{f2_ex0}	
\end{figure}

\noindent When $f_2=0$, there exists a unique constant-speed travelling front. In the singular limit, the speed is given by $c\approx 0.83$. Introducing $f_2$, we expect that outside the support of $f_2$ the dynamics are close to the dynamics of the constant-coefficient model; in particular, the singular limit of the front velocity will be close to $0.83$. We simulate the PDE for various values of $0<\varepsilon \ll 1$ (using MATLAB's \texttt{pdepe}), with an initial condition $(U_0,V_0)$ starting at position $z_0=-2.5$ (outside the support of $f_2$). That is, $z_0=-2.5$ is the starting position, in the sense that $U_0(z_0)=0$. The initial condition is constructed by first simulating the constant-coefficient PDE.
The fronts travel to the right with positive speed, and pass through the support of $f_2$. The front speed remains positive, but is not constant and fluctuates when the front is in the vicinity of the support of $f_2$.\\

\begin{figure}[h]
	   \centering 
       \scalebox{.8}{	\includegraphics[trim = 2cm 9.5cm 4cm 10cm, clip,width=0.8\textwidth]{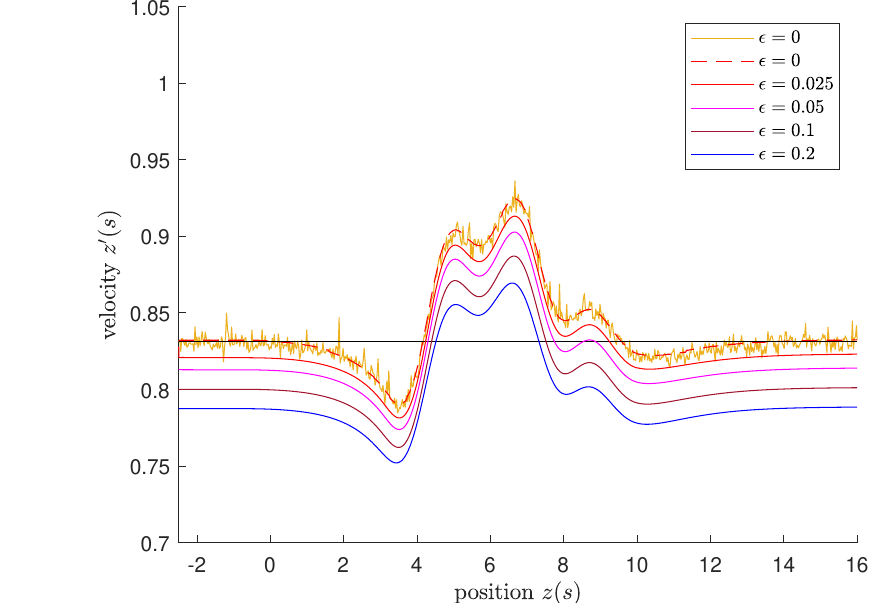}}
	\caption{ \small  Simulated position against velocity for the singular limit $\varepsilon=0$ (delay equation) and for various $\varepsilon>0$. For $\varepsilon=0$ the figure shows both the Monte Carlo approach (Section~\ref{s:delay_equation1}) to simulating the delay-differential equation, as well as the approach in Section~\ref{s:delay_equation2} based on the delay-differential equation in implicit form. The two approaches closely agree in the present example. 
    For simulations of the delay-differential equation using the Monte Carlo approach we use $M=10^5$ Monte Carlo samples, and a time discretization of $30$ steps per time unit $s=1$.   } 	
	\label{posvelocity_ex0}	
\end{figure}

\begin{figure}[H]
	   \centering 
          \scalebox{.8}{
       	\includegraphics[trim = 2cm 9.5cm 4cm 10cm, clip,width=0.8\textwidth]{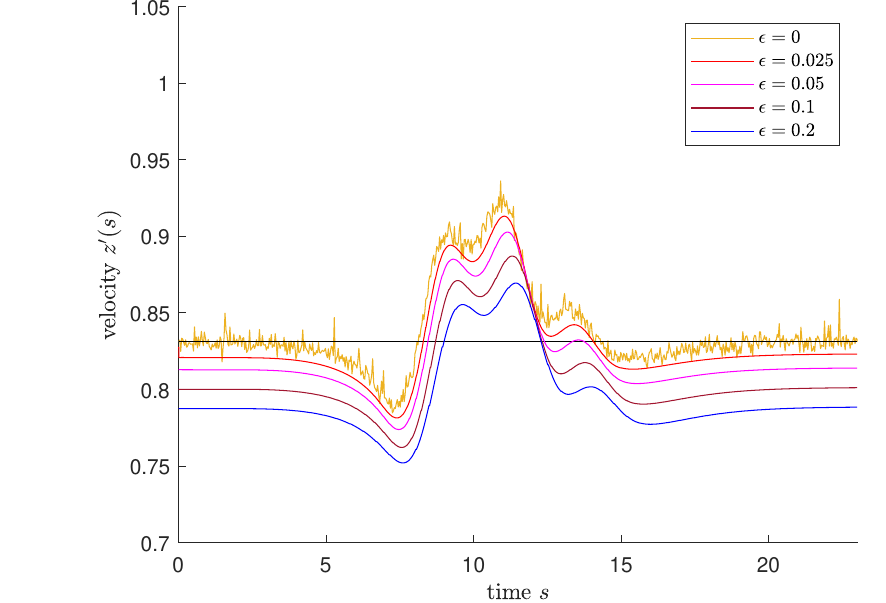}}
	\caption{ \small  Time against simulated velocity for the singular limit $\varepsilon=0$ (delay equation) and for various $\varepsilon>0$. }	
	\label{timevelocity_ex0}	
\end{figure} 

\paragraph{Example 1: Travelling front for $\alpha<0$ in the presence of stationary fronts.}
The aim of the present example is to illustrate Result~\ref{res:TF_withSF}, describing the interaction of travelling fronts and stationary fronts for positive feedback $\alpha<0$. To this end we consider the parameter values $\alpha=-2$, $\gamma=-0.2$, $\hat{\tau}=1$, and
\begin{align}
\begin{split} 
f_2(x)=& 0.3\exp(-0.5(x+0.1)^2)+0.33\exp(-2(x+1.5)^2)-0.53\exp(-2(x+0.75)^2)\\  
&+ 0.25\exp(-4(x-0.1)^2) - 0.4\exp(-3(x-1)^2)\, .
\end{split}
\end{align}

\begin{figure}[H]
	\centering 
	\scalebox{.9}{
	\includegraphics[trim = 2cm 9.5cm 4cm 10cm, clip,width=0.8\textwidth]{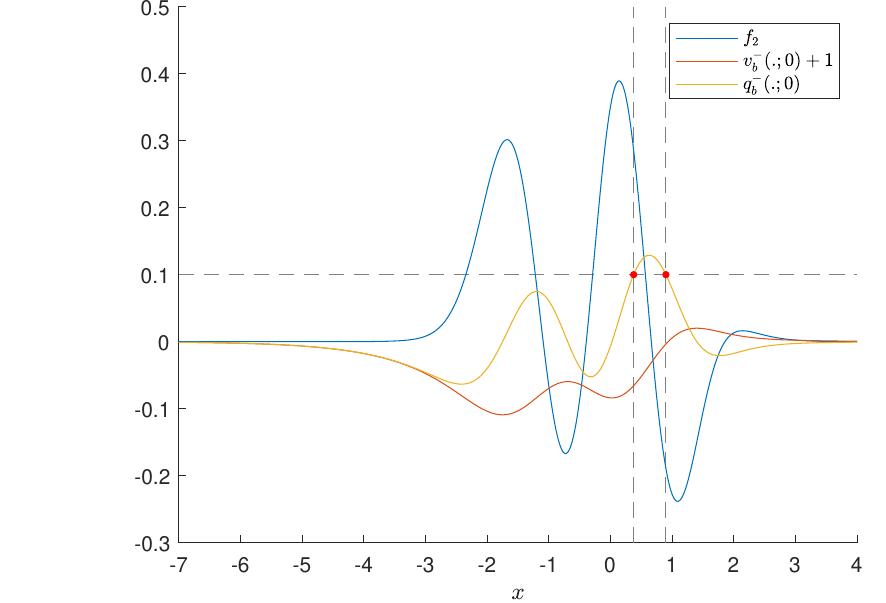}}
	\caption{ \small  $f_2(x)$, $v_b^-(x;0)+1$ and $q_b^-(x;0)$. The figure also indicates the singular limit positions of the stationary fronts given by the intersection points of $q_b^-(\cdot;0)$ with the horizontal line $\frac{\gamma}{\alpha}=0.1$. }	
	\label{f2_ex1}	
\end{figure} 
	
\noindent For $0<\varepsilon \ll 1$, there exist two stationary fronts $\textnormal{SF}_i(\varepsilon)$ with positions $z_{0,i}(\varepsilon)$, $i=1,2$, where by Result~\ref{res:stationaryfronts} the singular limit positions $z_{0,i}(0)$ are given by the solutions $z_0$ of the equation $q_b^-(z_0;0)=\frac{\gamma}{\alpha}=0.1$. Approximately, we have $z_{0,1}(0)\approx 0.38$ and $z_{0,2}(0)\approx 0.90$. The stationary front at position $z_{0,1}(\varepsilon)$ is unstable (since $(q_b^-)^\prime(z_{0,1}(0);\delta_0,0)>0$) and the stationary front at position $z_{0,2}(\varepsilon)$ is stable (see { Result~\ref{res:TF_withSF}}).\\

\noindent To generate an initial condition at a given target position $z_0 \in \mathbb{R}$ (in particular in the vicinity of the unstable stationary front for $\varepsilon>0$), we use the following algorithm. Let $z_0\in \mathbb{R}$ be the target position where we want to construct an initial condition $(U_0,V_0)$. Start with an initial condition $(U_{0,0},V_{0,0})$ which has a front-like shape (i.e.\ connecting the background states and having a single interface). Taking a sequence $(T_i)_i$, $T_i>0$, we iteratively construct initial conditions $(U_{0,i},V_{0,i})$ by solving the PDE~\eqref{eq:PDE_model_fast_reaction} with initial condition $(U_{0,i-1},V_{0,i-1})$ up to time $T_i$, and then shifting the resulting solution $(U,V)_{i-1}(x,T_i)$ in space to the position $z_0$, i.e.\ $(U_{0,i},V_{0,i})(x):= (U,V)_{i-1}(x-z_0+z_{i-1}(T_i),T_i)$. The main idea is to let $T_i\to 0$ while also ensuring that $\sum_i T_i$ is not too small in order to filter out transient dynamics. \\
	
\noindent We choose to work with the following sequence of times
$T_1=1$, $T_2=1$, $T_3=1$, $T_4=0.5$, $T_5=0.5$, $T_6=0.3$, $T_7=0.2$.\\

\noindent For $\varepsilon=0.05$ we estimate that the position $z_{0,1}(0.05)$ of the stationary front is contained in the interval $[0.37577,0.37609]$. The estimate arises from the fact that we constructed an initial condition $(U_{0,a},V_{0,a})$ at position $0.37577$ which ends up generating a front travelling to the left, and an initial condition $(U_{0,b},V_{0,b})$ at position $0.37609$ which generates a front travelling to the right. Likewise, for $\varepsilon=0.1$ we estimate that $z_{0,1}(0.1)$ is contained in the interval $[0.37714,0.38452]$.\\

\noindent The position of the stable stationary front is (surprisingly) a bit harder to estimate. On the one hand, stability implies that approaching fronts from the right and left are attracted, and consequently it appears that the position can be estimated in this way by tracking approaching fronts numerically for long enough. However, the numerics involve a discretization of space, and the position of the numerical front is confined to this grid. Hence, to estimate the position of the stationary front more precisely, one needs to increase the resolution of the grid, which quickly increases computation time.\\

\begin{figure}[!h]
	\centering
	\scalebox{.8}{%
		\includegraphics[trim = 1cm 8cm 2cm 8cm, clip,width=0.8\textwidth]
        {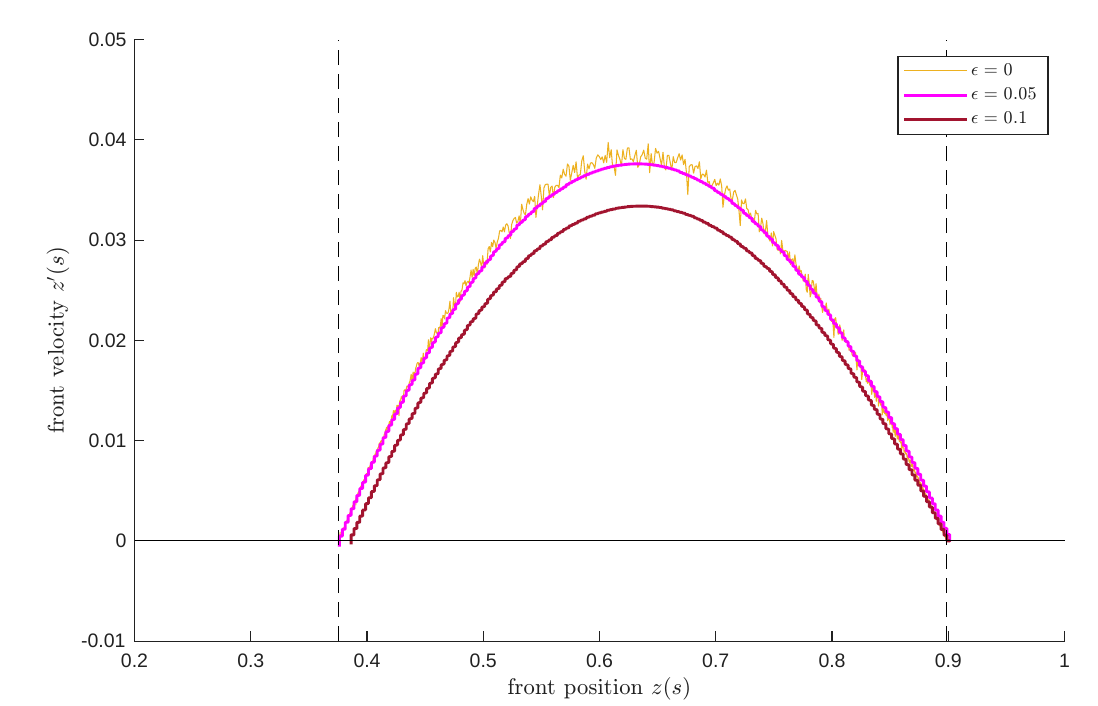}%
        \begin{picture}(0,0)
			\put(-300,100){ $SF_1$}
			\put(-40,100){ $SF_2$}
		\end{picture}%
	}
	\caption{\small Simulated position against velocity, delay-differential equation, $\varepsilon=0.05$ and $\varepsilon=0.1$. The { dashed} vertical lines indicate the singular limit positions of the stationary fronts. For simulation of the delay-differential equation we use Algorithm~1 in Section~\ref{s:delay_equation1} with $M=10^5$ Monte Carlo samples, and a time discretization of $30$ steps per time unit $s=1$.}
	\label{posvelocitya_ex1}
\end{figure}
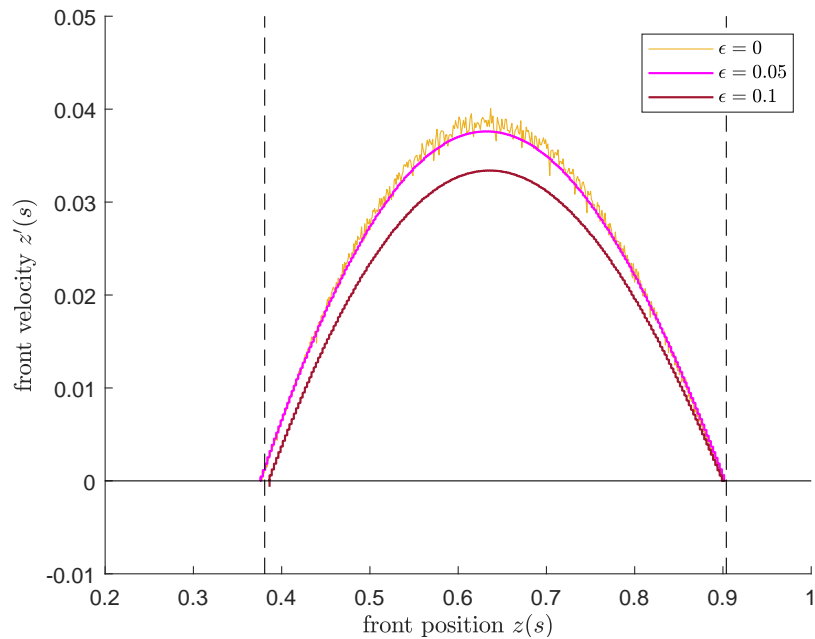

\begin{figure}[H]
	\centering
	\scalebox{.9}{\includegraphics[trim = 1cm 8.5cm 2.5cm 8cm, clip,width=0.8\textwidth]{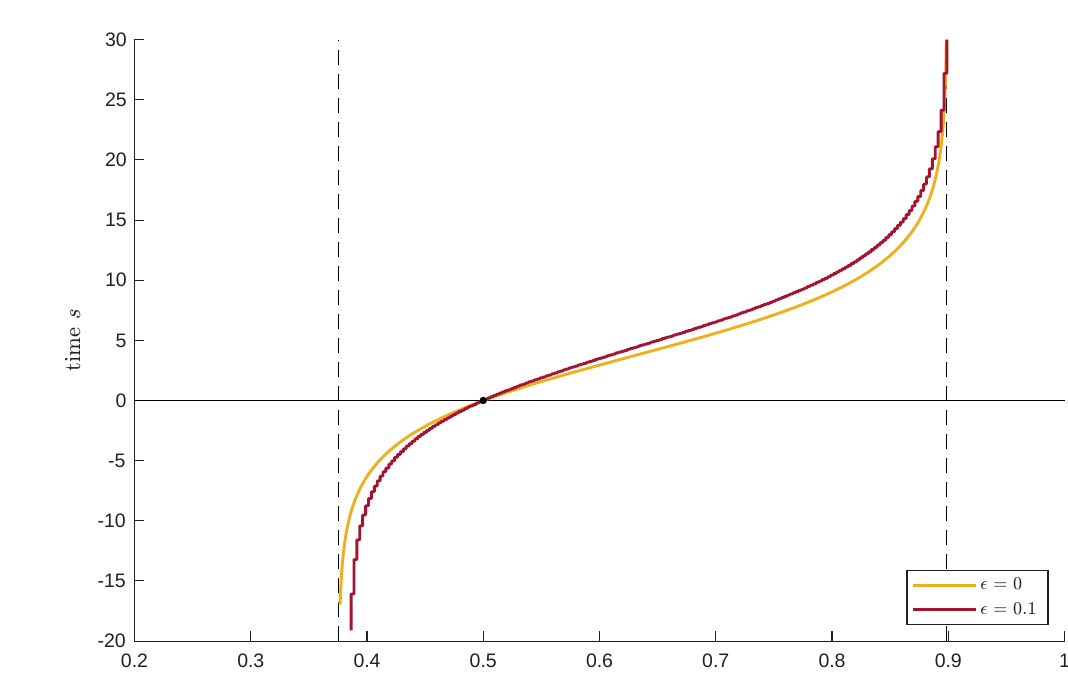}}
	\caption{ \small Simulated position against time, delay-differential equation and $\varepsilon=0.1$. The { dashed}  lines indicate the singular limit positions of the stationary fronts. To compare the delay-differential equation to the PDE~\eqref{eq:PDE_model_fast_reaction} for $\varepsilon=0.1$, we arrange (using translational symmetry in time) that the (time/position) trajectories pass through the common point $(z,s)=(0.5,0)$. }	
	\label{timepositiona_ex1}	
\end{figure}

\begin{figure}[H]
	\centering 
	\scalebox{.9}{\includegraphics[trim = 1cm 8.5cm 2.5cm 8cm, clip,width=0.8\textwidth]{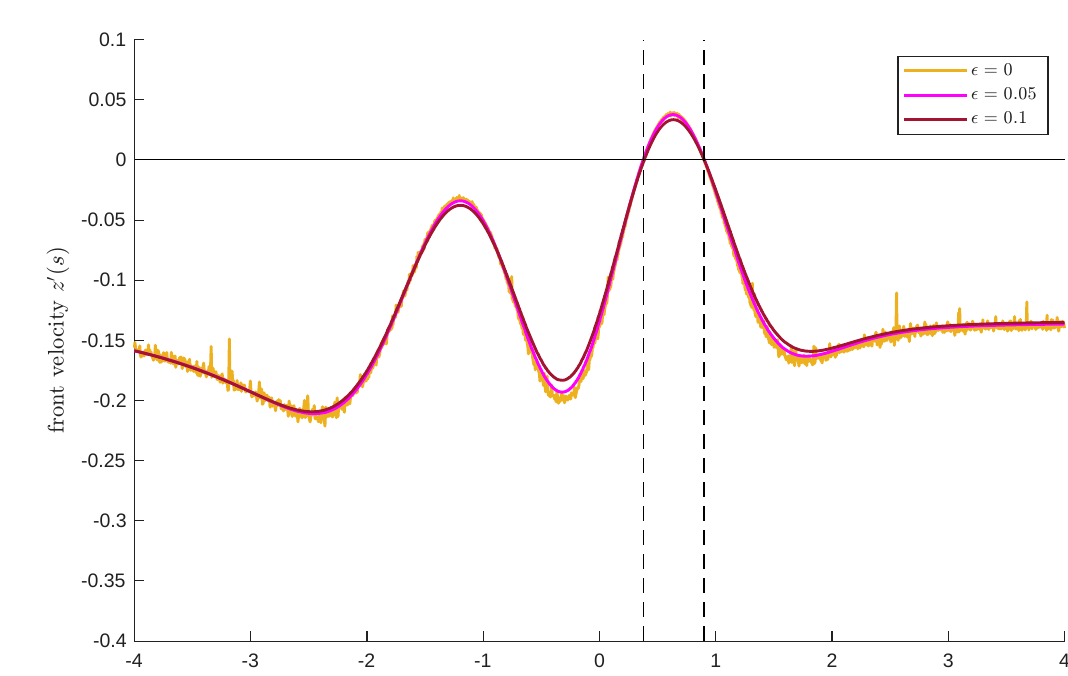}}
	\caption{ \small Simulated position against velocity, delay-differential equation, $\varepsilon=0.05$ and $\varepsilon=0.1$. 
	In addition to the travelling fronts connecting the unstable and stable stationary fronts shown in Figure~\ref{posvelocitya_ex1}, the present figure shows additional travelling fronts. Specifically, a travelling front which travels away from the unstable stationary front $\textnormal{SF}_1$ to the left, and a travelling front which approaches the stable stationary front $\textnormal{SF}_2$ from the right.}
	\label{posvelocity_ex1b}	
\end{figure}

\begin{figure}[H]
	\centering 	
	\scalebox{.9}{\includegraphics[trim = 1cm 8.5cm 2cm 7cm, clip,width=0.8\textwidth]{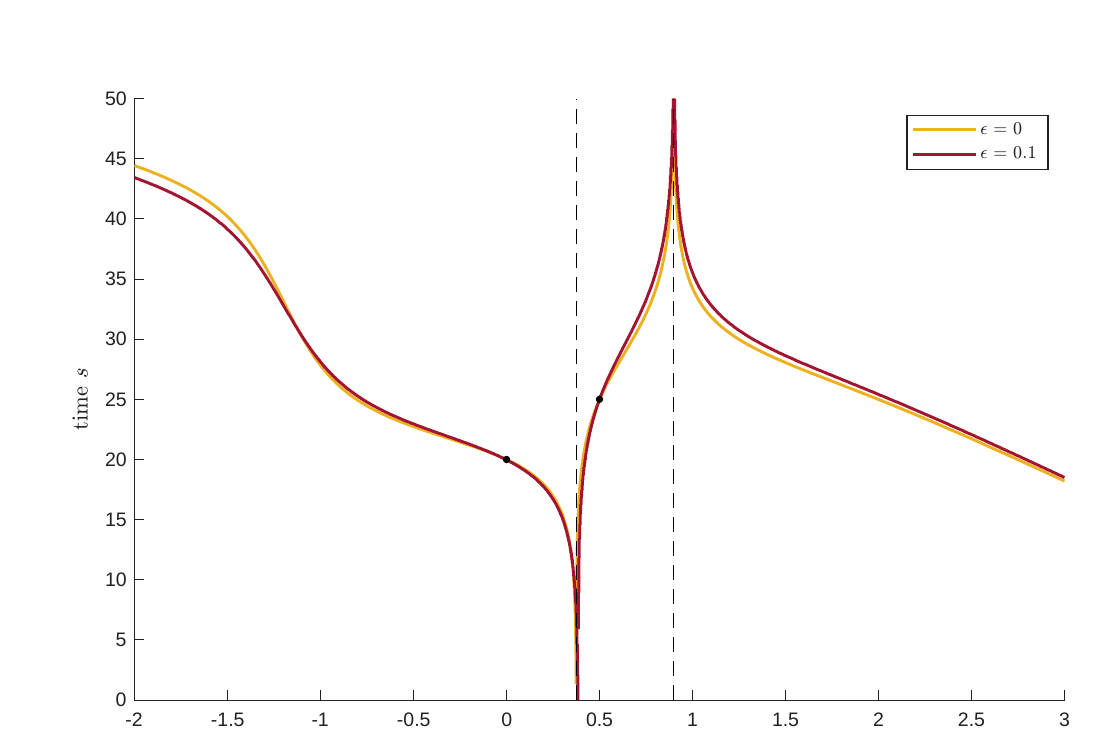}}
	\caption{ \small Simulated position against time, delay-differential equation, $\varepsilon=0.05$ and $\varepsilon=0.1$. To compare the delay-differential equation to the PDE~\eqref{eq:PDE_model_fast_reaction} for $\varepsilon=0.1$, we arrange (using translational symmetry in time) that the (time/position) trajectories pass through the common point $(z,s)=(0.5,25)$. }	
	\label{timeposition_ex1b}	
\end{figure}

\noindent Figure~\ref{posvelocitya_ex1} and Figure~\ref{timepositiona_ex1} illustrate a travelling front connecting stationary fronts $\textnormal{SF}_1$ in backward time and $\textnormal{SF}_2$ in forward time. Both the position $z(s)$ and the front velocity $z^\prime(s)$ closely agree with the associated singular position and velocity given by simulating the delay-differential equation.
In particular, when looking only at the front position $z(s)$ (Figure~\ref{timepositiona_ex1}), the comparison between $\varepsilon=0.1$ and $\varepsilon=0$ is already remarkably precise, which is why $\varepsilon=0.05$ is not included in this figure. On the other hand, the comparison for the front velocity $z^\prime(s)$ with the delay-differential equation (Figure~\ref{posvelocitya_ex1}) noticeably improves for $\varepsilon=0.05$ compared to $\varepsilon=0.1$.
Figures~\ref{posvelocity_ex1b} and Figure~\ref{timeposition_ex1b} extend the simulation by including other travelling fronts besides the travelling front connecting the stationary fronts. 
Specifically, a travelling front is included which travels away from the unstable stationary front $\textnormal{SF}_1$ to the left, and a travelling front which approaches the stable stationary front $\textnormal{SF}_2$ from the right.
    
\paragraph{Example 2: Direction-reversing front for $\alpha>0$ and $f_2=0$.}
    
\noindent We set $\hat{\tau}=1$, $\gamma=0.2$.  
Let $\alpha_{\textnormal{bp}}$ be given by 
\begin{align}
\alpha_{bp} =\left(\frac{\sqrt{2}}{3}c_{bp}- \gamma\right)\frac{ \sqrt{4+c_{bp}^2\hat{\tau}^2}}{c_{bp}\hat{\tau}} \approx 1.489\, , \quad 
c_{bp}= -\left(\frac{12\gamma}{\sqrt{2}\hat{\tau}^2}\right)^{\frac{1}{3}}\, .
\end{align}	
\noindent For $\alpha<\alpha_{bp}$, the algebraic equation \eqref{numerics_Melnikov} has a unique solution $c$ which corresponds to (up to translation in space) a unique constant-speed travelling front for the PDE~\eqref{eq:PDE_model_fast_reaction} for $0<\varepsilon \ll 1$. For $\alpha>\alpha_{bp}$, \eqref{numerics_Melnikov} has three solutions $c_m<c_0<c_p$, which for $0<\varepsilon \ll 1$ correspond to travelling fronts with speeds $c_m(\varepsilon)<c_0(\varepsilon)<0<c_p(\varepsilon)$.\\
	
\noindent For the present example, we take $\alpha=2.5>\alpha_{bp}$. We also fix $\varepsilon=0.1$ and compare it with the singular limit $\varepsilon=0$.\\
    
\noindent One can show that the delay-differential equation \cite[Chapter 5]{LvV_thesis} has an entire solution $z(s)$, which satisfies $z^\prime(s)\to c_0<0$ as $s\to -\infty$, $z^\prime(s)\to c_p>0$ as $s\to \infty$, and $z^\prime(s)$ is increasing in $s$. That is, in the singular limit $\varepsilon=0$ there exists a travelling front which changes direction, i.e.\ it initially travels to the left with negative speed (asymptotic to $c_0$ in backward time) and eventually approaches the positive speed $c_p>0$. In this example we explore numerically whether the PDE~\eqref{eq:PDE_model_fast_reaction} has a similar direction-reversing travelling front solution for $0<\varepsilon \ll 1$ (approaching $c_0(\varepsilon)$ and $c_p(\varepsilon)$ in backward, respectively forward time).\\
    
\noindent For setting up the simulation, the main challenge is to construct the initial condition $(U_0,V_0)$, which should be close to the front profile of the unstable constant-speed travelling front with speed $c_0(\varepsilon)$. To this end we consider the singularly perturbed system of ODEs for $(u,p,v,q)$ corresponding (in the fast scaling $\xi$) to the co-moving frame $\xi =\frac{x-cs}{\varepsilon}$. For a given value of $c$, there exists a unique orbit (standardized such that $u(0)=0$), which is forward asymptotic to the (independent of $x$, since in the present example $f_2=0$) background state $(u_b^+(\varepsilon),0,u_b^+(\varepsilon),0)$ and a unique orbit which is backward asymptotic to $(u_b^-(\varepsilon),0,u_b^-(\varepsilon),0)$. For $c=c_m(\varepsilon),c_0(\varepsilon),c_p(\varepsilon)$, these forward and backward orbits match at $z_0=0$, and form a heteroclinic. We numerically determine the value $c_0(\varepsilon)$ which gives rise to a heteroclinic orbit. To this end we identify locally invariant manifolds $\mathcal{M}_\varepsilon^\pm$ and restrict first to the one-parameter family of orbits which lie in the intersection of $W^u(\mathcal{M}_\varepsilon^-)$ and $W^s(\mathcal{M}_\varepsilon^+)$. To find this family of orbits we use an approximation of $\mathcal{M}_\varepsilon^\pm$ (the first order approximation in $\varepsilon$ already works well), and a `shooting' approach (tracking multiple orbits and investigating how closely they approach $\mathcal{M}_\varepsilon^\pm$). The resulting family is described by $u_0=0$ and $(p_0,v_0)$ being a function of $q$. We then use a second layer of `shooting' to investigate, for fixed $c$, which forward orbit in (the one-parameter family) $W^s(\mathcal{M}_\varepsilon^-)\cap W^u(\mathcal{M}_\varepsilon^+)$ converges to the positive background state, and which backward orbit converges to the negative background state. Lastly, using a bisection approach we can compute the value of $c$ where the forward and backward orbits match. According to the numerical results, $c_0(\varepsilon)$ (for $\varepsilon=0.1$) is contained in the interval $(-0.309375,-0.30875)$ (for comparison the singular limit velocity is given by $c_0(0)=-0.261792$). We now consider the following cases.\\
    
\noindent To construct an initial condition $(U_0,V_0)$ such that the solution of the PDE~\eqref{eq:PDE_model_fast_reaction} will reverse the direction of travelling, we take $c$ slightly larger than $c_0(\varepsilon)$ and define $(U_0,V_0)$ by means of concatenating the forward orbit converging to the positive background state with the backward orbit converging to the negative background state. How long the velocity $z^\prime(s)$ stays close to $c_0(\varepsilon)$ before moving away from it depends on how close $c$ is to $c_0(\varepsilon)$, as well as on other numerical factors (i.e.\ in the shooting algorithm described above). \\
 
\noindent To construct an initial condition $(U_0,V_0)$ for which the solution of the PDE will not reverse direction but instead converge to the stable travelling front with speed $c_m(\varepsilon)$, we take $c$ slightly smaller than $c_0(\varepsilon)$ and again concatenate forward and backward orbits.\\
   
\noindent For $\varepsilon=0$, we simulate the delay-differential equation, specifying the history to be $z_0=0$ and $z^\prime(s)=c_0(0)$ for $s\le 0$. Using the Monte Carlo approach for evaluating the delay functional, it depends on the specific sample realizations (i.e.\ seed for the pseudo random number generator) whether the velocity $z^\prime(s)$ ends up converging to $c_p(0)$ or $c_m(0)$. How long the solution of the delay-differential equation stays close to $c_m(0)$ before splitting up depends (in particular) on the number of Monte Carlo samples, and the time discretization.\\
    
\begin{figure}[H]
	\centering 
	\scalebox{.9}{
	\includegraphics[trim = 1cm 7.5cm 2cm 8cm, clip,width=0.8\textwidth]{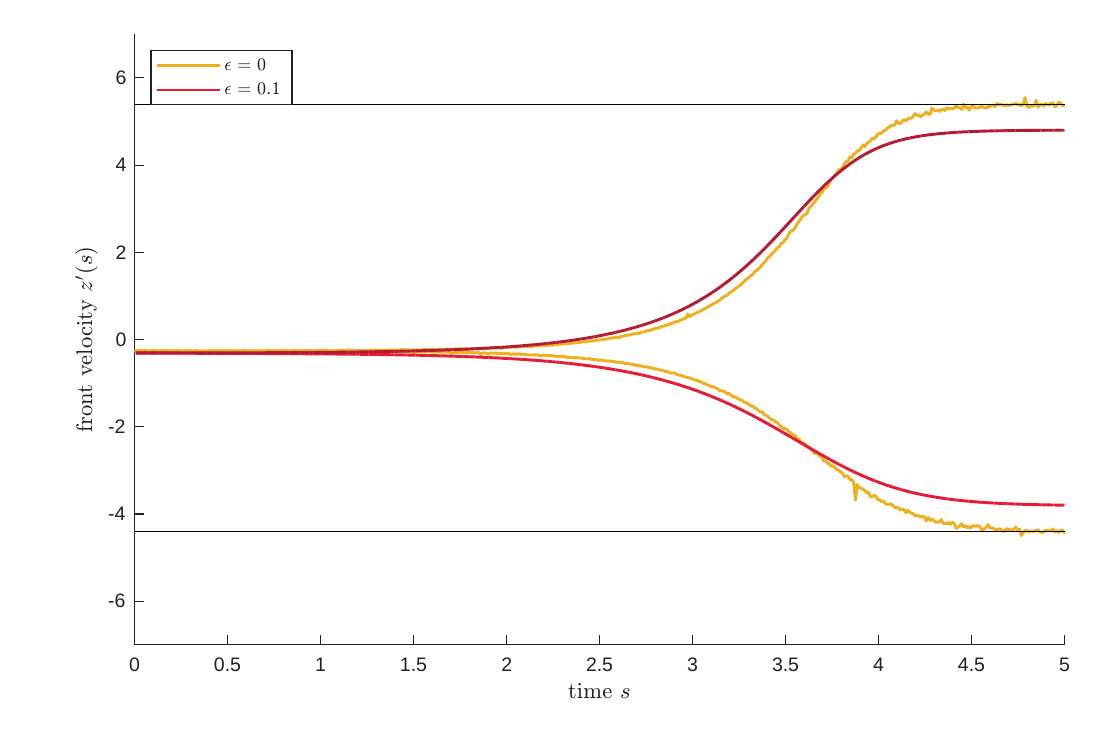}}
	\caption{ \small Simulated velocity against time, delay-differential equation and $\varepsilon=0.1$. For simulation of the delay-differential equation we use Algorithm~1 in
	Section~\ref{s:delay_equation1} with $M=1.5\times 10^5$ Monte Carlo samples, and a time discretization of $100$ steps per time unit $s=1$. The horizontal lines indicate the singular limit of the stable speeds $c_m(0)<0$ and $c_p(0)>0$. }	
	\label{timevelocity_ex2}	
\end{figure}

\begin{figure}[H]
	\centering
	\scalebox{.9}{\includegraphics[trim = 1cm 8.5cm 2cm 8cm, clip,width=0.8\textwidth]{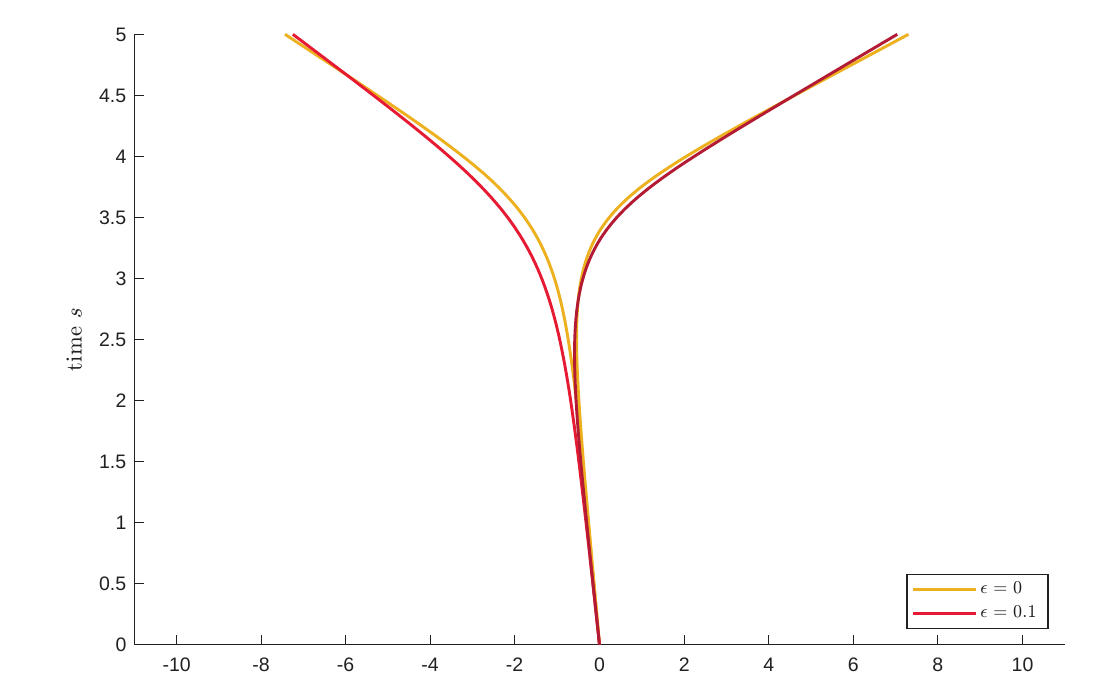}}
	\caption{ \small Simulated position against time, delay-differential equation and $\varepsilon=0.1$. }	
	\label{timeposition_ex2}	
\end{figure}

\begin{figure}[H]
	\centering 
	\scalebox{.9}{\includegraphics[trim = 1cm 7.5cm 2cm 7.5cm, clip,width=0.8\textwidth]{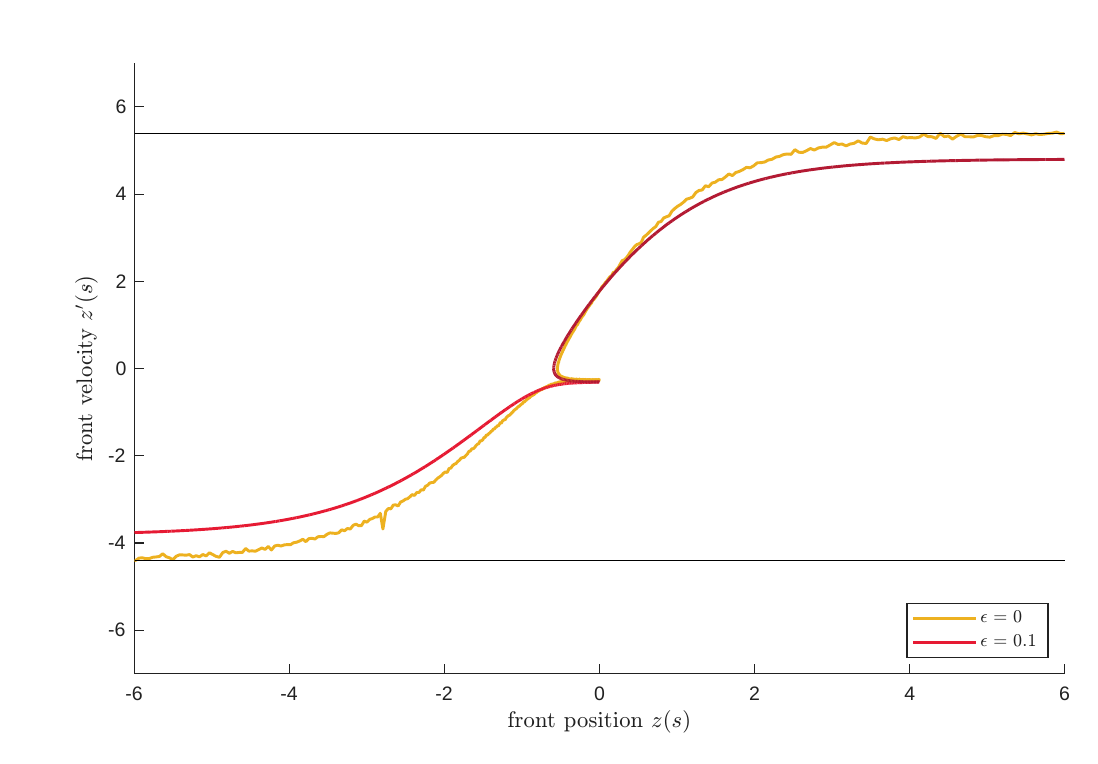}}
	\caption{ \small Simulated position against velocity, delay-differential equation and $\varepsilon=0.1$. The horizontal lines indicate the singular limit of the stable speeds $c_m(0)<0$ and $c_p(0)>0$. }	
	\label{posvelocity_ex2}	
\end{figure}

\noindent In summary, we obtain two numerical PDE solutions for $\varepsilon=0.1$, one which connects the velocities $c_0(\varepsilon)<0$ and $c_p(\varepsilon)>0$ (and in particular reverses its direction of travelling) and a second solution which connects the velocities $c_0(\varepsilon)<0$ and $c_m(\varepsilon)<0$. For $\varepsilon=0$ we obtain similar solutions of the delay-differential equation. Figures~\ref{timevelocity_ex2}, \ref{timeposition_ex2}, \ref{posvelocity_ex2} depict the position and velocity compared for $\varepsilon=0$ and $\varepsilon=0.1$.

\paragraph{Example 3: Direction switches and periodic front.}
\noindent A periodic front for the parameter values $(\alpha,\gamma,\hat{\tau})$ of Example 2 when the inhomogeneity $f_2$ is introduced.\\

\noindent Unlike Examples 0 and 1, which illustrate theoretical results (i.e.\ Result~\ref{res:delay_eq} and Result~\ref{res:TF_withSF}, respectively), and Example 2 which illustrates partial results (i.e.\ the existence of direction-reversing entire solutions for the delay-differential equation ($\varepsilon=0$), but not for the PDE \eqref{eq:PDE_model_fast_reaction} for $\varepsilon>0$), the present example is purely exploratory (in the sense that it is not backed up by theoretical results). In fact the $f_2$ which we will work with does not satisfy Assumptions~\ref{ass:f}.\\

\noindent The delay-differential equation is of the form
\[
\gamma + \alpha\bigl(\hat{\tau}W[z](s)-q_b^-(z(s);0)\bigr)=\frac{\sqrt{2}}{3}z^\prime(s).
\]
In particular, it appears that sudden large changes in $q_b^-$ might have a similar impact on the velocity $z^\prime(s)$. This raises the question whether a localised inhomogeneity $f_2$ might act as a `direction switch', i.e.\ cause approaching fronts to reverse their direction.\\
	
\noindent We explore this question in the context of the parameters in Example 2, where we discussed that in the absence of $f_2$ there exist stable speeds $c_m(\varepsilon)<0$ and $c_p(\varepsilon)>0$. We introduce the following $f_2$ given by two localized pulses at positions $z=\pm 11$ 
\begin{align} 
f_2(x)= - \rho \exp(-40(x+11)^2) - \rho\exp(-40(x-11)^2)\, .
\end{align} 
		
\noindent We set $\rho=12$. The graph of the corresponding $q_b^-(\cdot;0)$ has large changes in the vicinity of positions $z=\pm 11$ of the pulses and is close to $0$ everywhere else.

\begin{figure}[H]
	\centering
	\scalebox{.9}{\includegraphics[trim = 1cm 7.5cm 2cm 8cm, clip,width=0.8\textwidth]{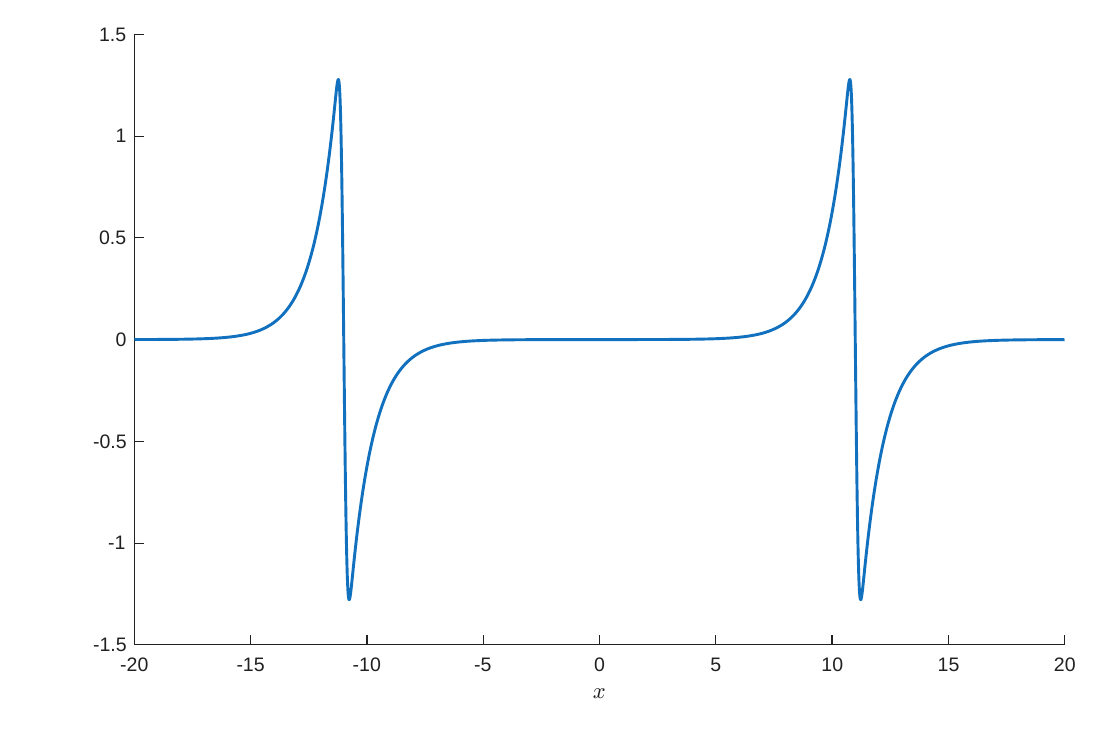}}
	\caption{ \small Graph of $q_b^-(x;0)$. The graph exhibits two localized regions where it changes in a rapid way. The intention of the design is to construct direction switches which reverse the direction of approaching fronts, and trap them in a compact space given by the space between the two switches. The choice for the scaling $\rho=12$ in the definition of $f_2$ comes from the fact that it makes the switches strong enough for this purpose. For example, reversing a front approaching the Gaussian pulse at $z=-11$ from the left seems to require taking $\rho$ around $4$.
        }	
	\label{qmbin_ex3}
\end{figure} 	

\noindent Near $z_0=0$ (far away from $\pm 11$) we expect the front dynamics to be close to the constant-coefficient model where $f_2=0$. Consider a front travelling to the left away from $z_0=0$ with speed close to $c_m(\varepsilon)<0$. This front approaches $z=-11$ from the right and will encounter the large negative local minimum of $q_b^-$. Given the sudden decrease of $q_b^-$, $z^\prime(s)$ might quickly increase as a response, and if the drop in $q_b^-$ is large enough, $z^\prime(s)$ might become positive, i.e.\ the front reverses direction. \\
	
\noindent Likewise, consider a front travelling to the right away from $z_0=0$ and towards $z=11$ with speed close to $c_p(\varepsilon)$. This front will encounter the large local maximum of $q_b^-$. Consequently, the front speed $z^\prime(s)$ might decrease, and could become negative.\\

\begin{figure}[H]
	\centering 	
	\scalebox{.9}{\includegraphics[trim = 1cm 7.5cm 2cm 8cm, clip,width=0.8\textwidth]{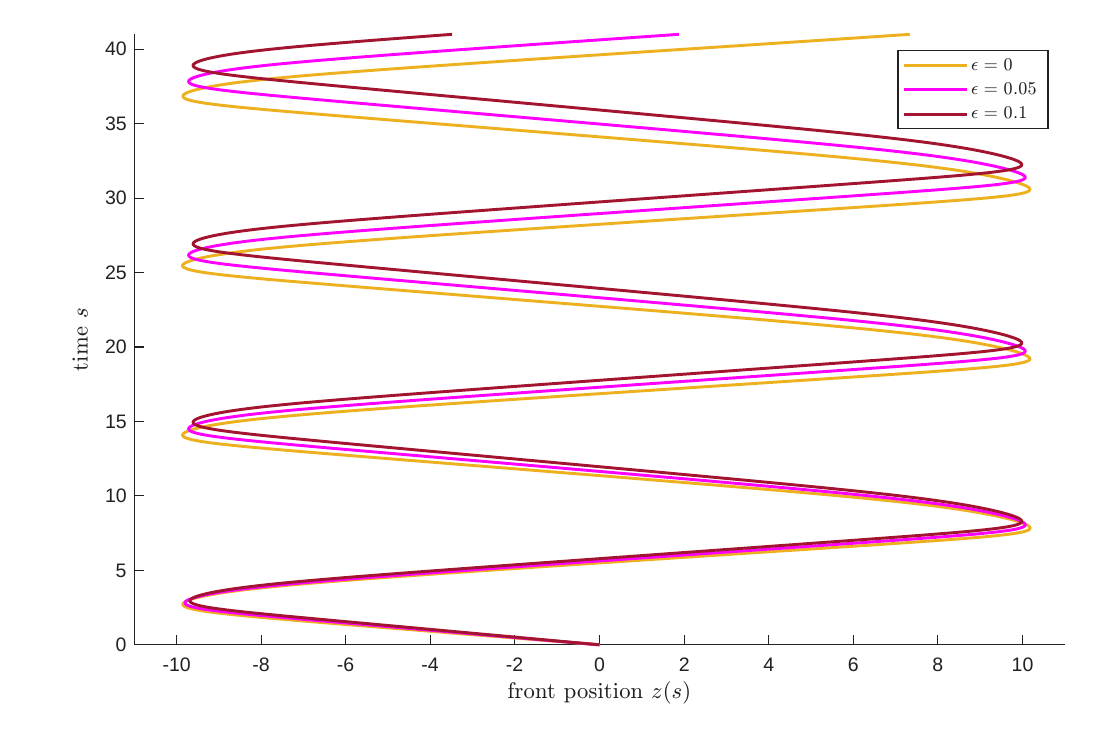}}
	\caption{ \small Simulated position against time, delay-differential equation, $\varepsilon=0.05$ and $\varepsilon=0.1$. }	
	\label{timeposition_ex3}	
\end{figure}

\noindent The numerical results confirm expectations. Figure~\ref{timeposition_ex3} shows the simulated front position $z(s)$ against time and suggests that the pulses placed at positions $z=\pm 11$ act as direction switches as intended. The change of direction occurs around $z=\pm 10$. In particular, the front position is trapped in the compact interval $[-11,11]$. A periodic front seems to emerge as a consequence of the two direction switches.
	
\begin{figure}[!h]
	\centering 	
	\scalebox{.9}{\includegraphics[trim = 1cm 7.5cm 2cm 8cm, clip,width=0.8\textwidth]{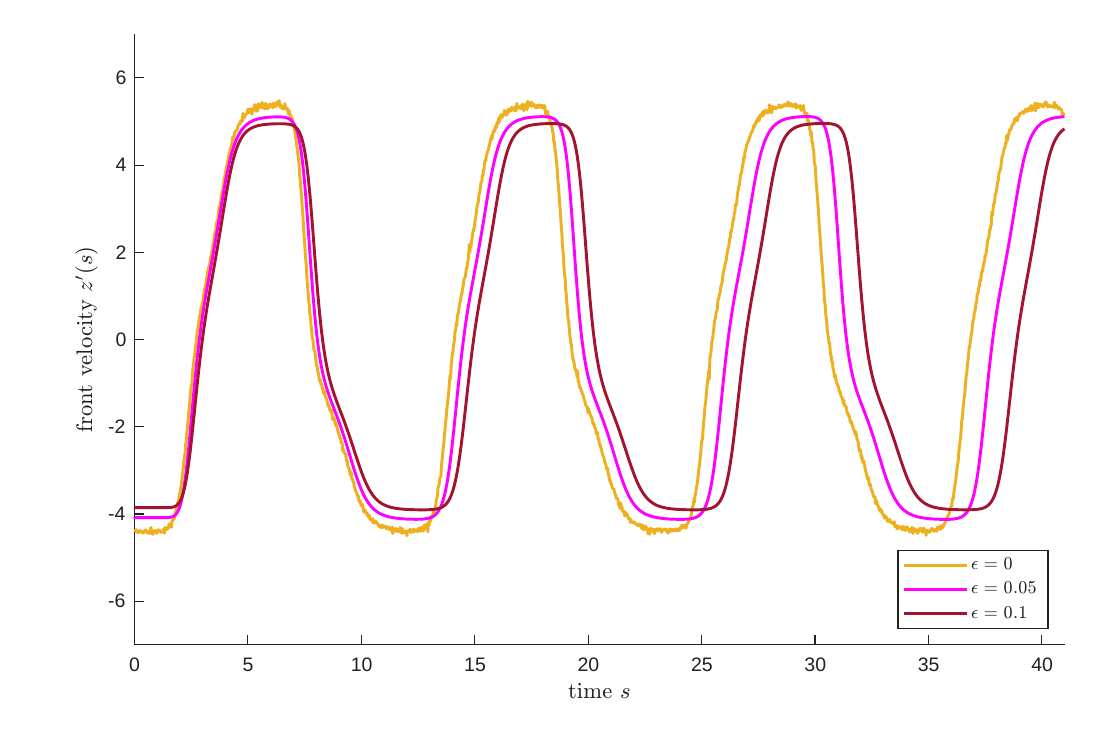}}
	\caption{ \small Simulated velocity against time, delay-differential equation, $\varepsilon=0.05$ and $\varepsilon=0.1$. For the delay-differential equation we use Algorithm~1 in
	Section~\ref{s:delay_equation1} with $M=1.5\times 10^5$ Monte Carlo samples and a time discretization of $40$ steps per time unit $s=1$. In addition we also explored Algorithm~2 in
	Section~\ref{s:delay_equation2} (not shown) based on the implicit form of the delay-differential equation, but unlike in Example 0, Algorithm~2 does not match the results of Algorithm~1 and in fact does not produce an orbit which appears periodic. The underlying theoretical and/or numerical reasons why Algorithm~2 exhibits behaviour that appears `irregular' are currently unknown. 
	}	
	\label{timevelocity_ex3}	
\end{figure} 
		
\begin{figure}[H]
	\centering 
	\scalebox{.9}{\includegraphics[trim = 1cm 7.5cm 2cm 8cm, clip,width=0.8\textwidth]{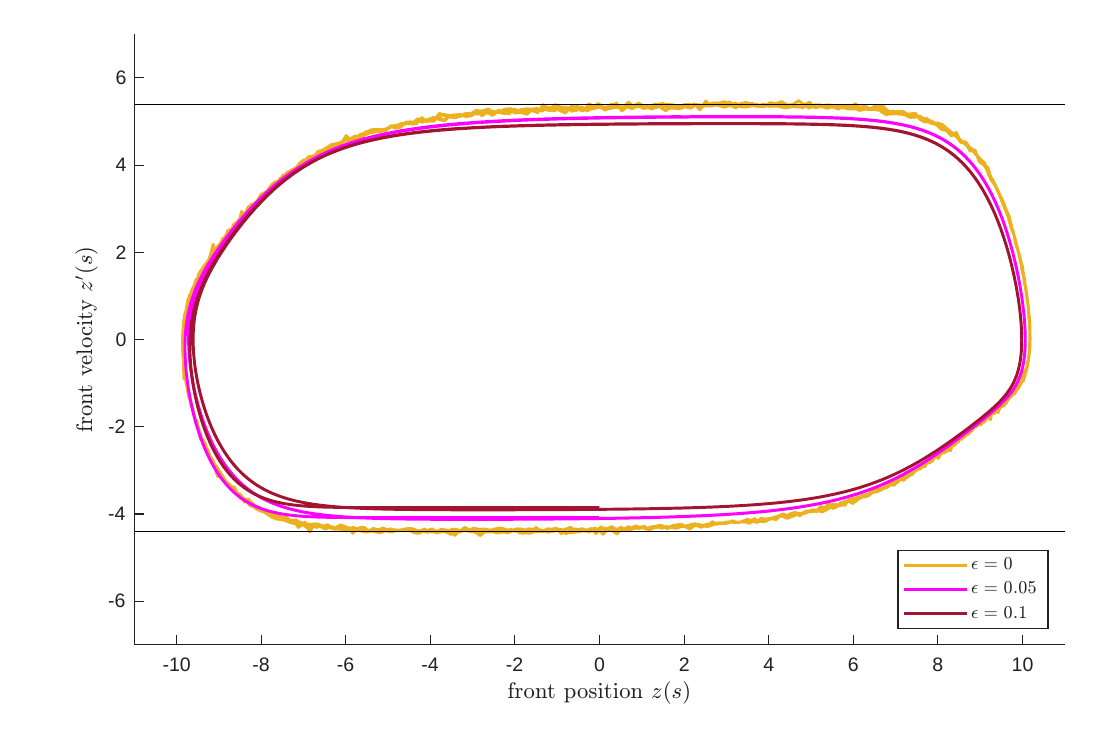}} 
	\caption{ \small Simulated position against velocity, delay-differential equation, $\varepsilon=0.05$ and $\varepsilon=0.1$. }	
	\label{posvelocity_ex3}	
\end{figure} 

\begin{figure}[H]
	\centering 
	\scalebox{.9}{\includegraphics[trim = 1cm 7.5cm 2cm 8cm, clip,width=0.8\textwidth]{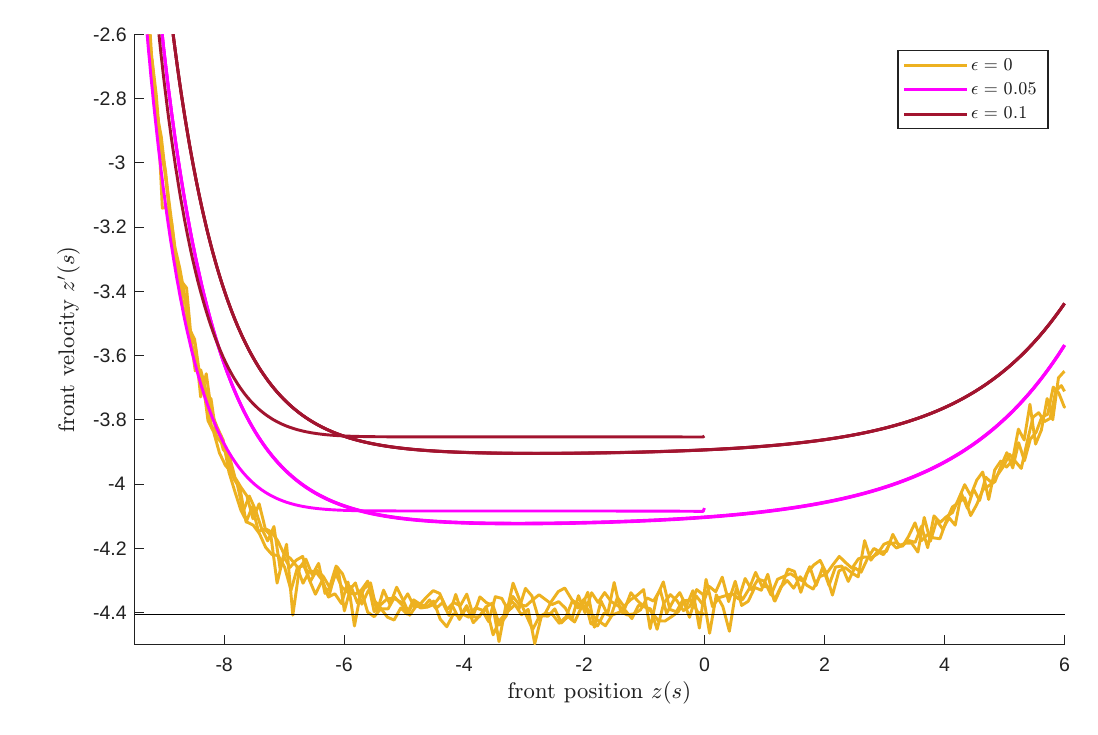}}
	\caption{ \small Simulated position against velocity, delay-differential equation, $\varepsilon=0.05$ and $\varepsilon=0.1$. The figure shows a zoomed-in version of Figure~\ref{posvelocity_ex3}. The figure shows that (due to transient dynamics) the initial part of the front position/velocity trajectory for $\varepsilon=0$ and $\varepsilon=0.1$ deviates from the periodic orbit which seems to be approached eventually. More specifically, the deviation is particularly noticeable when the position/velocity trajectory approaches the southwest region of the domain for the first time. Interestingly, the simulated path for the delay-differential equation ($\varepsilon=0$) does not appear to have a noticeable initial transient phase. }	
	\label{posvelocity_ex3zoom}	
\end{figure}

\begin{figure}[H]
\centering
\scalebox{1}{
\begin{tikzpicture}
  [
  ->,
  >=stealth',
  auto,node distance=1cm,
  thick,
  main node/.style={circle, draw, font=\sffamily\Large\bfseries}
  ]
  
  \node[anchor=south west,inner sep=0] (image) at (0,0)
  { \includegraphics[trim = 2cm 9cm 4cm 9cm, clip,width=0.8\textwidth]{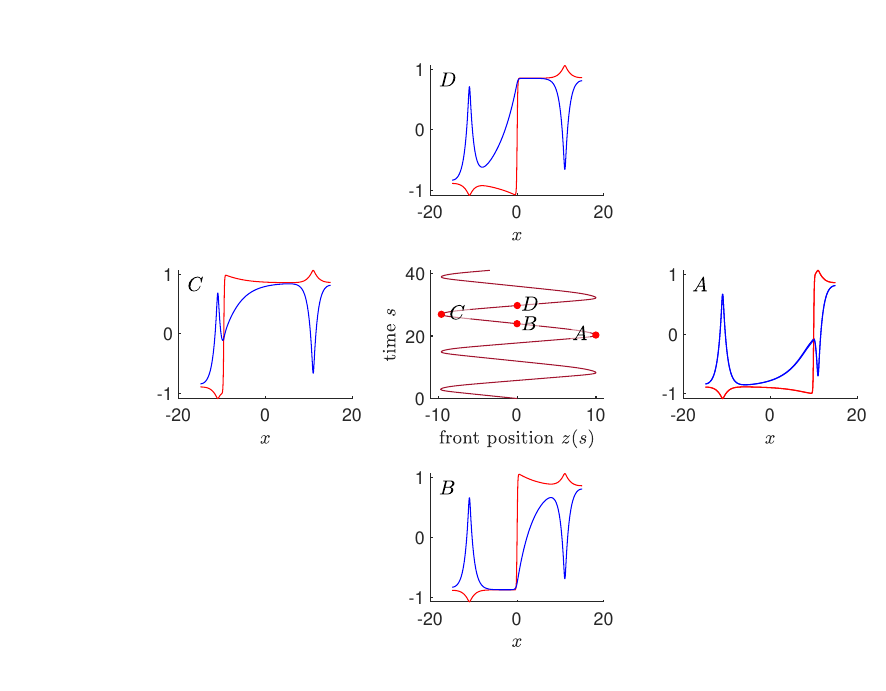} };
  
  \begin{scope}[x={(image.south east)},y={(image.north west)}]
    \node[align=center] (A1) at (0.9,0.3) { };
    \node[align=center] (A2) at (0.75,0.15) { };
    \draw [->] (A1) to [out=270,in=0] (A2);
    
    \node[align=center] (B1) at (0.4,0.15) { };
    \node[align=center] (B2) at (0.25,0.3) { };
    \draw [->] (B1) to [out=180,in=270] (B2);	
    
    \node[align=center] (C1) at (0.25,0.65) { };
    \node[align=center] (C2) at (0.4,0.8) { };
    \draw [->] (C1) to [out=90,in=180] (C2);					
  \end{scope} 
\end{tikzpicture}
}
\caption{ \small Front profile of solution $(U,V)$ at various points $A$, $B$, $C$, $D$ of the simulated periodic orbit for $\varepsilon=0.1$. The $U$-component is shown in red and the $V$-component is shown in black. }	
\label{frontprofile_ex2}	
\end{figure}
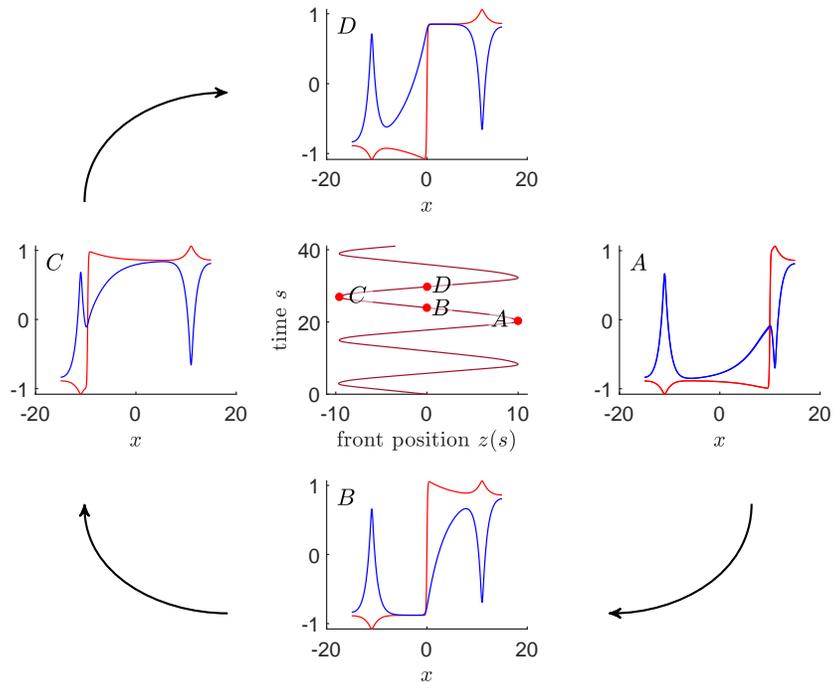 

\newpage


\section{Summary and outlook}
We have demonstrated that \eqref{eq:PDE_model} exhibits bi-stability, characterized by two stable stationary background states that vary in space due to the presence of $f_1, f_2$ and further give rise to stationary and travelling front solutions that connect them.

\medskip

The construction of stationary background states in the stationary wave ODE \eqref{eq:background_ODEsystem} needed an extension of Fenichel theory to the general non-autonomous case (that is, going beyond \cite{Eldering.2012} to include also unstable directions) combined with the theory of exponential dichotomies. Using a Melnikov-type argument for \eqref{eq:background_ODEsystem} gave also a leading order existence condition for stationary front solutions \eqref{gammaSF0} that resembles the one for the constant-coefficient case \eqref{eq:existence_const_coeff}, the difference being that, for fixed parameter settings, the possible front positions $x_0$ are fixed (contrasting with the translation-invariant case of constant-coefficient equations). Travelling front solutions have a richer structure then their counterparts in the constant-coefficient case: their interfaces move non-uniformly through the motionless heterogeneous background states, which makes the usual travelling wave ansatz hopeless. Upon carefully defining them along their stationary version and numerical observations (see Definition~\ref{def:front}), it turns out that one can prove (i) their existence from a specific set-up of the corresponding initial value problem (leading to initialised fronts as in Definition~\ref{def:Initialised_front}) and (ii) their existence for all $t \in \mathbb{R}$ (leading to entire front solutions as in Definition~\ref{def:entire_front}). The main novelty of the proof for the IVP lies in tracking the time evolution of the time-derivatives $(\partial_t U, \partial_t V)$ which are bi-asymptotic to zero as $x \to \pm \infty$ for each fixed $t$ (since the background states are motionless even for travelling front solutions). The proof of the existence of entire solutions is built on the result for the IVP and uses the construction of a localised co-moving frame and the Arzela-Ascoli theorem. The exposition closes with a short formal derivation of a delay-differential equation \eqref{eq:delay_eq} that describes the temporal evolution of the front position (in the sense of Definition~\ref{def:position}). We briefly sketch a numerical algorithm to solve it and give various examples, both for cases where its validity can be proven rigorously and cases where its validity is expected but not proven. All proofs can be found in \cite{LvV_thesis} and will also be subject of follow-up articles for which the present article serves as overview.

\medskip

Future research will further generalize the obtained results to other settings with particular interest in treating also nonlinear slow equations, that is, the $+U$ in the $V$-equation in \eqref{eq:PDE_model} replaced by $F(U,V)$. Furthermore, a careful study of the dynamics exhibited by the delay-differential equation \eqref{eq:delay_eq} will be of interest, but is expected to be challenging due to the difficult type of delay differential equation. Another immediate question is its relation with the traditional approach of deriving reduced ODEs for the dynamic position using a rigorous center manifold reduction (CMR)(for a "collective coordinates type" approach). This is especially interesting since the derivation of \eqref{eq:delay_eq} made no smallness assumption on $z$, which one does when constructing travelling fronts as bifurcating from stationary ones. Note that, the conventional CMR approach for \eqref{eq:PDE_model} would need control over the spectrum of the operator obtained from linearisation around stationary fronts. This is, in general, difficult unless one enforces smallness assumptions for $f_1, f_2$ or assumes them to be periodic or spatially localised (as done in \cite{BC-BD.2020} for the stability of pulses). The derivation of the delay-differential equation completely circumvents this difficulty.

\medskip

Further capitalising on the derivation of the delay-differential equation that avoids intricate spectral problems, a promising direction is to use this approach for interacting multi-front and pulse solutions under the influence of spatial heterogeneity. Going beyond the 1-D case, a promising direction is the study of 2-D front solutions and more complicated 2-D patterns such as the tracking of spiral wave tips in reaction-diffusion systems (as given, e.g. in \cite{Barkley1994PRL, BarkleyKevrekidis1994Chaos,SandstedeScheelWulff1997,Wulff1998ZPC}) with applications, e.g. in the modelling of cardiac dynamics where heterogeneities are believed to play a major role. 

\newpage

\bibliographystyle{arXiv_article_v2}
\bibliography{arXiv_article_v2}

\end{document}